%% file: agt-5-42.tex
\documentclass{gtart_h}
\input agtout

\lognumber{42}
\volumenumber{5}
\volumeyear{2005}
\papernumber{42}
\pagenumbers{1027}{1050}
\received{21 September 2004} 
\accepted{2 August 2005}
\published{25 August 2005}

\usepackage{amsmath,amssymb,graphicx,psfrag}

\def\psfraga <#1,#2> #3#4{%
\psfrag {#3}{\smash{\rlap{\kern #1 \raise #2\hbox{#4}}}}}

\def\fref#1{\hyperlink{#1anchor}{#1}}
\def\figref#1{\hyperlink{#1anchor}{Figure~#1}}
\def\anchor#1{\noindent\hypertarget{#1anchor}{\smash{$\phantom{99}$}}}

\usepackage{}
\theoremstyle{plain}
\newtheorem{theorem}{Theorem}

\newtheorem{proposition}[theorem]{Proposition}
\theoremstyle{remark}
\newtheorem{remark}{Remark}

\begin{document}
\title{A search method for thin positions of links}
\authors{Daniel J. Heath\\Tsuyoshi Kobayashi}

\address{Mathematics Department,  Pacific Lutheran University\\Tacoma, WA
98447, USA\\{\rm and}\\Department of Mathematics, Nara Women's 
University\\Kitauoya Nishimachi, Nara 630, Japan} 
\asciiaddress{Mathematics Department,  Pacific Lutheran University\\Tacoma, WA
98447, USA\\and\\Department of Mathematics, Nara Women's 
University\\Kitauoya Nishimachi, Nara 630, Japan}

\asciiemail{ heathdj@plu.edu, tsuyoshi@cc.nara-wu.ac.jp }
\gtemail{\mailto{heathdj@plu.edu},
\mailto{tsuyoshi@cc.nara-wu.ac.jp}}

\begin{abstract}
We give a method for searching for thin 
positions of a given link.
\end{abstract}


\maketitle

\section{Introduction}

The concept of {\it thin position} of a knot or link was introduced by
David Gabai (see \cite{G}) in 1987, and has since played an important
role in 3--dimensional topology (see for example \cite{G-L},
\cite{S-T}).  However, it seems to be quite difficult to find a thin
position for a given link.  In \cite{T}, Abigail Thompson shows that
either a link in thin position is also in bridge position, or the link
has an essential meridional planar surface properly imbedded in its
complement.  This allows us to determine, up to the determination of
the bridge positions of the link, the thin positions of the link,
provided that the exterior of the link does not admit an essential
meridional planar surface.  This result is generalized by Yo'av Rieck
and Eric Sedgwick \cite{R-S} to show that any thin position of a
connected sum of small knots is obtained by placing minimal bridge
presentations of the factor knots vertically and taking their
connected sum by using monotonic arcs.  Moreover Ying-Qing Wu shows
that a thinnest level sphere in a thin position of a link gives an
essential meridional planar surface \cite{Wu}.  The purpose of this
paper is, by using the authors' previous result \cite{H-K}, to give a
search method for finding thin positions of a given link, up to the
determination of the bridge positions of given signed graphs and
meridional essential planar surfaces, which gives a natural
generalization of the result of Rieck-Sedgwick.

The first author is supported by a PLU Regency Advancement Award 
and a Wang Center Award. 
The second author is supported by a Grant-in-Aid for Scientific Research, Japan Society for the Promotion of Science.

\section{Preliminaries}

We begin with some definitions.

For a submanifold $H$ of a manifold $K$, $N(H,K)$ denotes a regular neighborhood of $H$ in $K$.

Let $L$ be a non-splittable link in $S^3$.
The {\it complement} of $L$, denoted by $E(L)$, is the closure of the complement of a regular neighborhood of $L$ in $S^3$. 
A {\it meridional} planar surface in the complement of $L$ is a planar surface properly imbedded in the link complement with boundary components consisting only of meridians of $L$.

Note that $S^3 \backslash \lbrace \text{two~points} \rbrace = S^2 \times {\mathbb R}$.
We define $p: S^2 \times {\mathbb R} \rightarrow S^2$ to be the projection onto the first factor, and $h:S^2 \times {\mathbb R} \rightarrow {\mathbb R}$ to be the projection onto the second factor.
We recall the definition of thin position. 
Suppose that $h \vert _L$ is a Morse function on $L$. 
Let $c_1, c_2, \dots , c_n$ $(c_1 < c_2 < \dots < c_n)$ be the critical values of $h \vert _L$. 
Choose regular value $r_i$ $(1 \le i \le n-1)$ so that $c_i < r_i < c_{i+1}$. 
Then the {\it width} of $L$, denoted by $w(L)$, is defined as follows. 
$$w(L) = \sum_{i=1}^{n-1} \vert L \cap h^{-1} (r_i) \vert$$
We say that $L$ is in a {\it thin position}, 
if for any ambient isotopy $f_t$ $(0 \le t \le 1)$ with 
$h \vert _{f_1(L)}$ a Morse function on $f_1(L)$, we have $w(L) \le w(f_1(L))$. 
We say that $S$ is a {\it thin 2--sphere (thick 2--sphere)} for $L$ with respect to $h$ if $S = h^{-1}(t)$ for some $t$ which lies between adjacent critical values $c_i$ and $c_{i+1}$ of $h \vert_{L}$, where $c_i$ is a maximum (minimum) of $L$, and $c_{i+1}$ is a minimum (maximum) of $L$. 
A 2--sphere $S$ in $S^2 \times {\mathbb R}$ is said to be {\it bowl like} if $S = F_1 \cup F_2$ such that $F_1 \cap F_2 = \partial F_1 = \partial F_2$, $F_1$ is a round 2--disk contained in a level 2-sphere, $h\mid_{F_2}$ is a Morse function with exactly one maximum or minimum, $p(F_1) = p(F_2)$, and $p\mid_{F_2}:F_2 \rightarrow p(F_2)$ is a homeomorphism.
Further, when we consider a link $L$ together with the bowl like 2-sphere, we shall require that all punctures by $L$ lie in $F_1$.
A bowl like 2--sphere is {\it flat face up (flat face down)} if $F_1$ is above (below) $F_2$ with respect to $h$.

A spacial graph $G$ is a 1-complex embedded in the 3-sphere. 
In particular, $G$ is a {\it signed vertex graph} if each vertex of $G$ is labelled with either $+$ or $-$.
We define the {\it width} of a signed vertex graph $G$. 
When we consider width of $G$, we suppose that 1) the vertices of $G$ labelled with $+$ ($-$ resp.) have the same height and are higher (lower resp.) than any other points in $G$, and 2) $h \vert_{G \setminus \{\text{vertices}\}}$ is a Morse function. 
We say that $G$ is in a {\it bridge position} if each maximum in $G \setminus \{\text{vertices}\}$ is higher than any minimum of $G \setminus \{\text{vertices}\}$.
In general, let $r_1, \dots , r_{n-1}$ $(r_1 < \dots < r_{n-1})$ be regular values between the critical values in $G \setminus \{\text{vertices}\}$. 
Then we define the width of $G$ by the following. 
$$w(G) = \sum_{i=1}^{n-1} \vert G \cap h^{-1}(r_i) \vert$$
We say that $G$ is in a {\it thin position} if for any ambient isotopy $f_t$ $(0 \le t \le 1)$ such that 

\begin{enumerate}
\item
All the points of $f_1( \text{the vertices labelled with $+$ ($-$ resp.)})$ have the same height and are higher (lower resp.) than any other points in $f_1(G)$, and 

\item
$h \vert_{f_1(G \setminus \{\text{vertices}\})}$ is a Morse function, 
\end{enumerate}

\noindent 
we have that $w(G) \le w(f_1(G))$.

For a signed vertex graph $G$ in a bridge position, we define the bridge number of the bridge position by  $\vert F \cap G \vert/2$, where $F$ is a level 2-sphere such that every maximum (minimum resp.) of $G$ is above (below resp.) $F$. 
The minimum of bridge numbers for all possible bridge positions of the signed vertex graph $G$ is the {\it bridge index} of $G$. 
Obviously, this is a generalization of the concept of bridge index for links.

In general, let $L$ be an unsplittable link, and $\mathcal{S} = S_1 \cup \dots \cup S_m$ a union of mutually disjoint bowl like 2-spheres with the following property (Property~1). 
Let $C_0, C_1, \dots, C_m$ be the closures of the components of $S^3 \setminus \mathcal{S}$ such that $C_0$ lies exterior to all of $S_j$, and $C_i$ $(i=1, \dots , m)$ is the component lying directly inside of $S_i$. 
Then for each $j$ $(j=0, 1, \dots, m)$, we suppose the following is satisfied (see Property~1 in page 109 of \cite{H-K}).

\noindent
{\bf Property 1}

\begin{enumerate}
\item 
For each $C_j$ $(j=0, 1, \dots , m)$, we have either one of the following. 

\begin{enumerate}
\item 
There are both a maximum and a minimum of $L$ in $C_j$, or 
\item
There does not exist a critical point of $L$ in $C_j$. 
\end{enumerate}

\item 
There exists a level 2-sphere $F_0$ in $C_0$ such that: 

\begin{enumerate}
\item  
every flat face down (up resp.) bowl like 2-sphere in $\partial C_0$ lies above (below resp.) $F_0$, and 
\item 
every maximum (minimum resp.) of $L$ in $C_0$ (if one exists) lies above (below resp.) $F_0$, and it is lower (higher resp.) than the flat face down (up resp.) bowl like 2-spheres in $\partial C_0$.
\end{enumerate}

\item 
For each $i$ $(i=1, \dots, m)$, there exists a level disk $F_i$ properly embedded in $C_i$ such that: 

\begin{enumerate}
\item 
every flat face down (up resp.) bowl like 2-sphere in $\partial C_i \setminus S_i$ lies above (below resp.) $F_i$, and 
\item 
every maximum (minimum resp.) of $L$ in $C_i$ (if one exists) lies above (below resp.) $F_i$, and it is lower (higher resp.) than the flat face down (up resp.) bowl like 2-spheres in $\partial C_i \setminus S_i$. 
\end{enumerate}
See \figref{2.1}. 
\end{enumerate}

\begin{figure}[ht]\small\anchor{2.1}
\begin{center}
\includegraphics[width=6cm, clip]{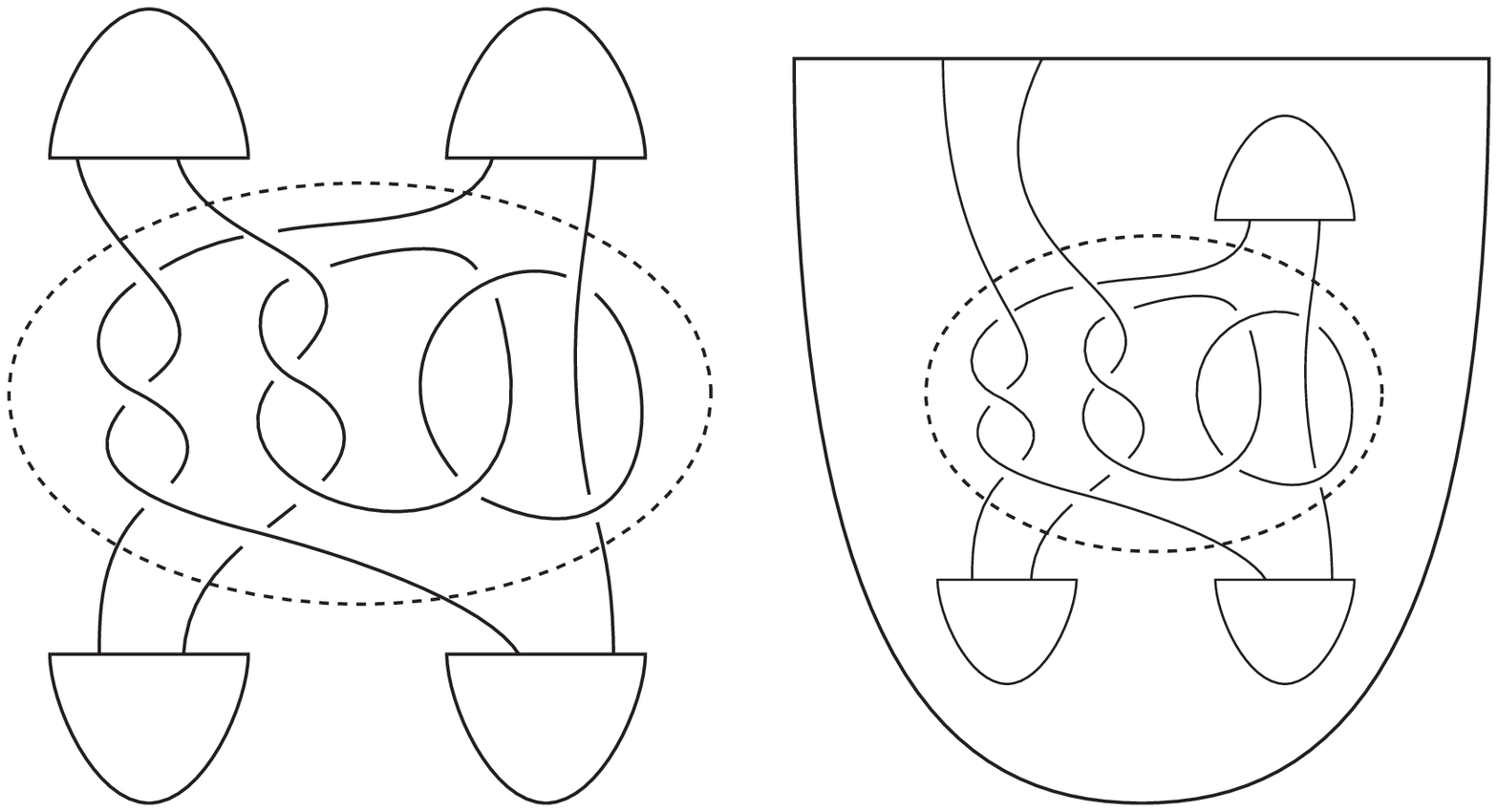}
\end{center}
\begin{center}
Figure 2.1
\end{center}
\end{figure}

\noindent
{\bf Signed vertex graph obtained from $(C_j, L \cap C_j)$}

Let $L$, $C_j$ $(j=0, 1, \dots , m)$ be as above. 
We can obtain a signed vertex graph $G_j$ from $(C_j, L \cap C_j)$ as follows. 

Suppose that $j=0$. 
In this case, by shrinking each component of $\partial C_0$ to a vertex, and by pulling up (down resp.) the vertices obtained from flat face down (up resp.) 2-spheres to make them in the same level, we obtain a signed vertex graph, say $G_0$, where the high (low resp.) vertices are the $+$ ($-$ resp.) vertices. 
By (2) of Property~1, we see that $G_0$ is in a bridge position, see \figref{2.2}.

\begin{figure}[ht]\small\anchor{2.2}
\begin{center}
\includegraphics[width=6cm, clip]{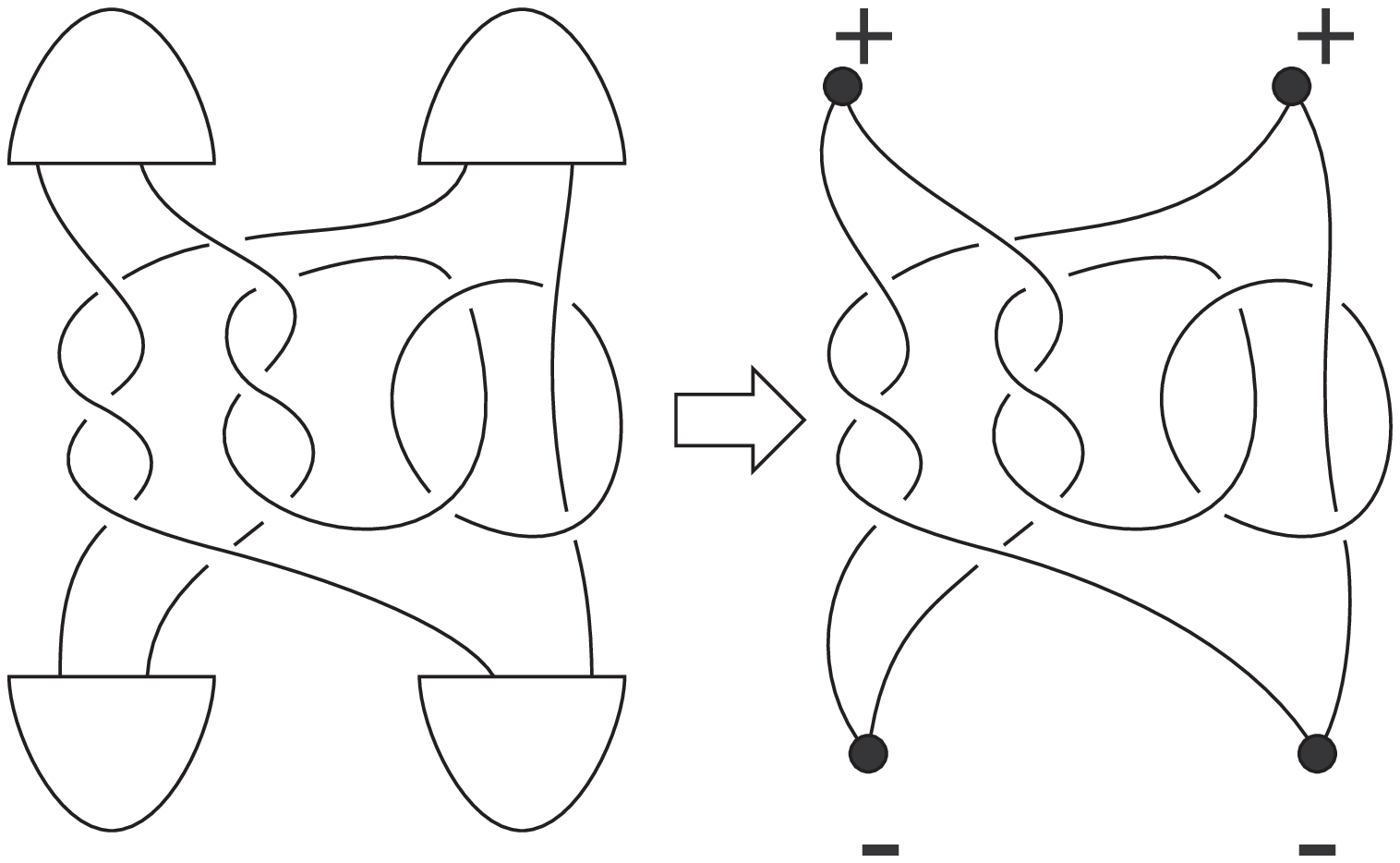}
\end{center}
\begin{center}
Figure 2.2
\end{center}
\end{figure}

\noindent
Suppose that $j \ne 0$. 
In this case we may deform $C_j$ by an ambient isotopy, $f_t$, of $S^3$ which does not alter $(\text{flat face of }S_j) \cup (L \cap C_j) \cup (\partial C_j \setminus S_j)$, making $f_1(C_j)$ appear as of type $C_0$, see \figref{2.3}. 
Notice that $f_t$ causes $S_j$ to pass through the ``points at infinity" of $S^3$ with respect to $h$. 
Then we apply the above argument for $\Big(f_1(C_j), f_1(L \cap C_j)\Big)$ to obtain a signed vertex graph $G_i$ in a bridge position. 
We say that $G_j$ $(j=0, 1, \dots , m)$ is a {\it signed vertex graph associated to} $\mathcal{S}$.

\begin{figure}[ht]\small\anchor{2.3}
\begin{center}
\includegraphics[width=6cm, clip]{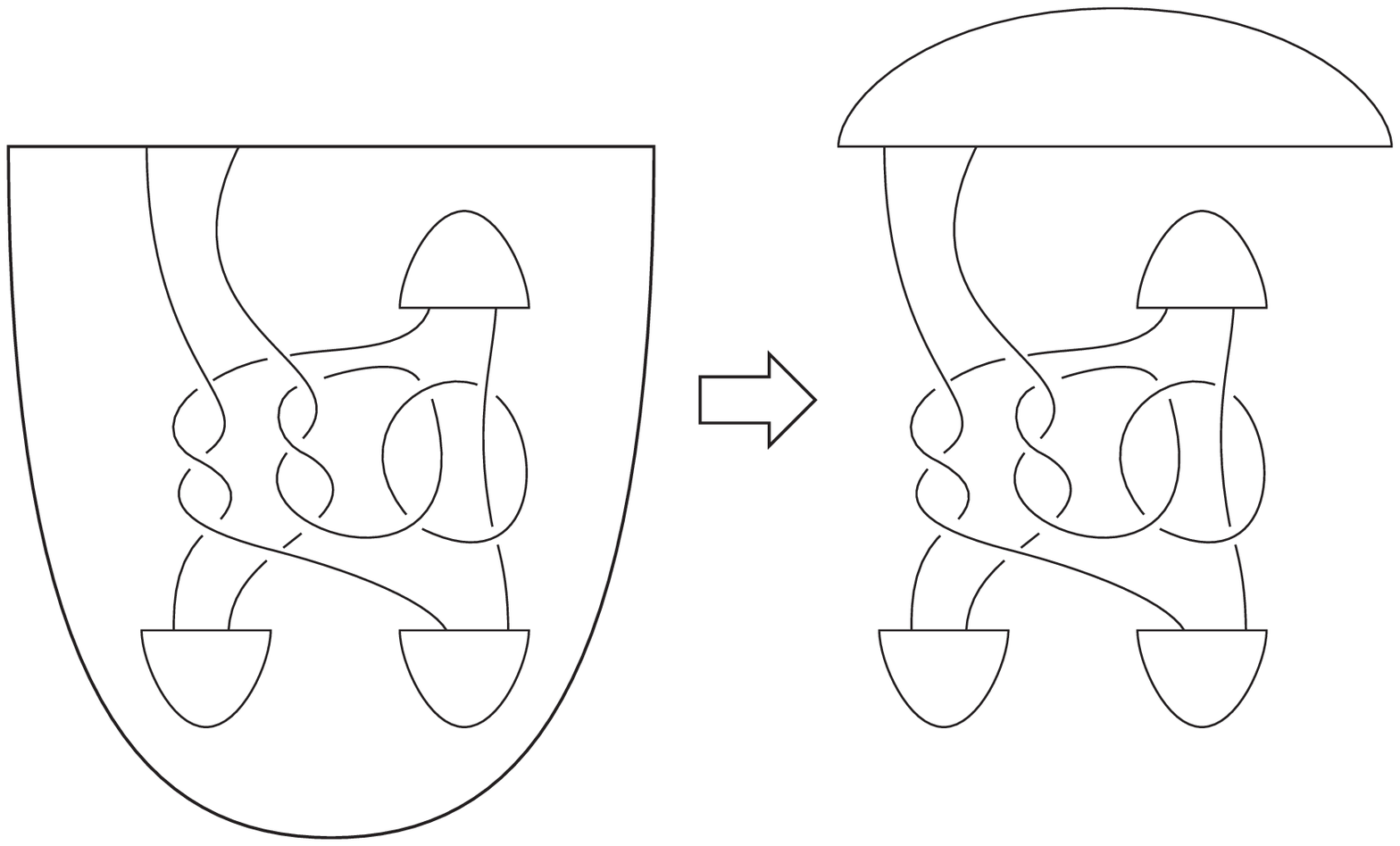}
\end{center}
\begin{center}
Figure 2.3
\end{center}
\end{figure}

\noindent
{\bf Cocoons}

Let $L$, $C_j$ $(j=0, 1, \dots , m)$ be as above. 
Then we can take a convex 3-ball $R_j$ in $\text{Int} C_j$ such that each component of $\text{cl} ((L \cap C_j)\setminus R_j)$ is a monotonic arc connecting $R_j$ and a component of $\partial C_j$, 
and that 

\begin{enumerate}
\item 
$R_0$ lies below (above resp.) the flat face down (up resp.) bowl like 2-spheres in $\partial C_0$. 

\item
$R_i$ $(i=1, \dots, m)$ lies below (above resp.) the flat face down (up resp.) bowl like 2-spheres in $\partial C_i \setminus S_i$.

\end{enumerate}

See \figref{2.1}. 
We call $R_j$ a {\it cocoon} of $L$ associated to $\mathcal{S}$. 

\begin{remark}
By the definition of the cocoons, we have the following. 

\begin{enumerate}

\item 
Each critical point of $L$ in $C_j$ (if one exists) is contained in $R_j$. 

\item 
Each component of $L \setminus \cup_{j=0}^m R_i$ is a monotonic arc intersecting $\mathcal{S}$ in exactly one point. 
\end{enumerate}

\end{remark}

Let $\prec$ be a linear order on the cocoons $\{ R_0, R_1, \dots , R_m\}$. 
We say that $\prec$ is {\it compatible with relative positons in} $L$ (associated to $\mathcal{S}$) if the following is satisfied. 

\begin{quote}
For each pair $i, j$ $(i \ne j)$ with $C_i \cap C_j \ne \emptyset$, $R_i \prec R_j$ if and only if $R_j$ is located above $R_i$ in $L$. 
\end{quote}

\noindent
{\bf Note}\qua 
By (2) of Remark~1, we see that if $C_i \cap C_j \ne \emptyset$, then $h(R_i)$ and $h(R_j)$ are disjoint. 

\begin{proposition}
Let $\{ R_0, R_1, \dots , R_m\}$, $\prec$ be as above. 
Suppose that $\prec$ is compatible with relative positons in $L$. 
Then there is an ambient isotopy $g_t$ $(0 \le t \le 1)$ which satisfies the following. 

\begin{enumerate}
\item 
For each $j$ $(j=0,1, \dots , m)$, $L \cap R_j$ and $g_1( L \cap R_j )$ are similar. 

\item 
Each component of $g_1(L\setminus \cup_{j=0}^m R_j)$ is a monotonic arc. 

\item 
The intervals $h(g_1( L \cap R_0 ))$, $h(g_1( L \cap R_1 ))$, \dots, $h(g_1( L \cap R_m ))$ are mutually disjoint, and the order of the positions of the intervals in the real line ${\mathbb R}$ agrees with $\prec$. 

\end{enumerate}
\end{proposition}

\begin{proof}
We first fix mutually disjoint $m + 1$ intervals $I_0, I_1, \dots , I_m$ whose order of the positions of the intervals in the real line agrees with $\prec$, i.e., $I_j$ is higher than $I_i$ if and only if $R_i \prec R_j$. 

\begin{figure}[ht]\small\anchor{2.4}
\begin{center}
\psfrag{B}{$B_{s_1}$}
\psfraga <-2pt,0pt> {I}{$I_0$}
\psfraga <-3pt,0pt> {R}{$R_0$}
\includegraphics[width=6cm, clip]{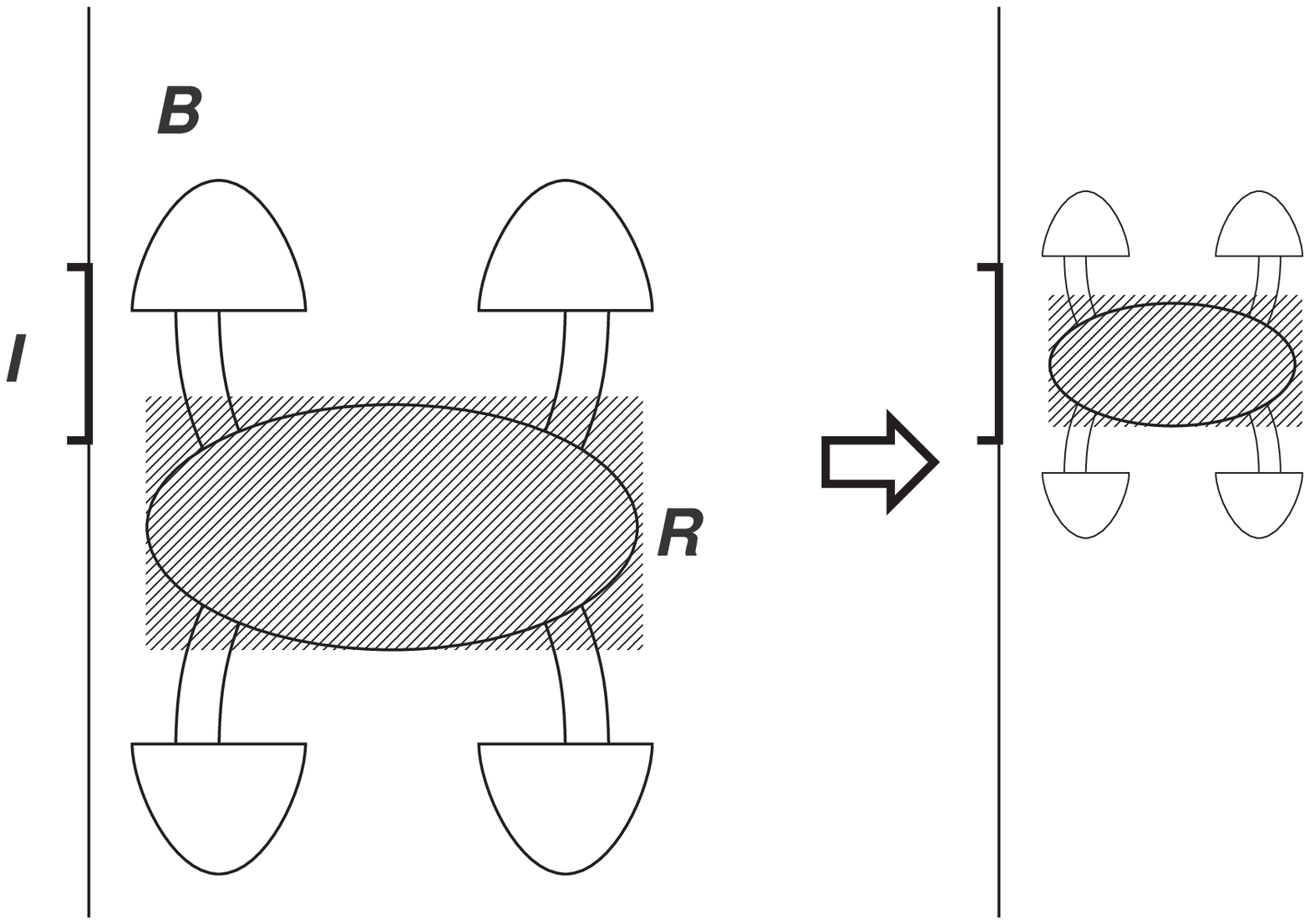}
\end{center}
\begin{center}
Figure 2.4
\end{center}
\end{figure}

We construct a desired ambient isotopy as follows. 
First, we deform $L$ by a similarity deformation and parallel translation so that $h( R_0 )$ is contained in $I_0$, see \figref{2.4}. 
For simplicity, we abuse notation by denoting the image of $L$ under this isotopy also by $L$. 
Let $S_{s_1}, \dots , S_{s_p}$ be the components of $\partial C_0$, and $B_{s_1}, \dots , B_{s_p}$ mutually disjoint 3-balls bounded by $S_{s_1}, \dots , S_{s_p}$ respectively. 
Then we next apply an ambient isotopy such that $B_{s_1}, \dots , B_{s_p}$ are deformed by a similarity transformation, so that the length of $h( B_{s_i} )$ is shorter than the length of $I_{s_i}$ $(1 \le i \le p)$. 
Then we pull up or down $B_{s_i}$ according to whether $S_{s_i}$ is flat face down or up, to make $h( B_{s_i} ) \subset \text{Int} I_{s_i}$, see \figref{2.5}. 
We note that since $\prec$ is compatible with relative positions in $L$, it is easy to see that we can make the ambient isotopy to satisfy the following. 

\begin{figure}[ht]\small\anchor{2.5}
\begin{center}
\psfrag{B1}{$B_{s_1}$}
\psfrag{B2}{$B_{s_2}$}
\psfrag{B3}{$B_{s_3}$}
\psfrag{B4}{$B_{s_4}$}
\psfrag{I1}{$I_{s_1}$}
\psfrag{I2}{$I_{s_2}$}
\psfrag{I4}{$I_{s_3}$}
\psfrag{I3}{$I_{s_4}$}
\includegraphics[width=6.5cm, clip]{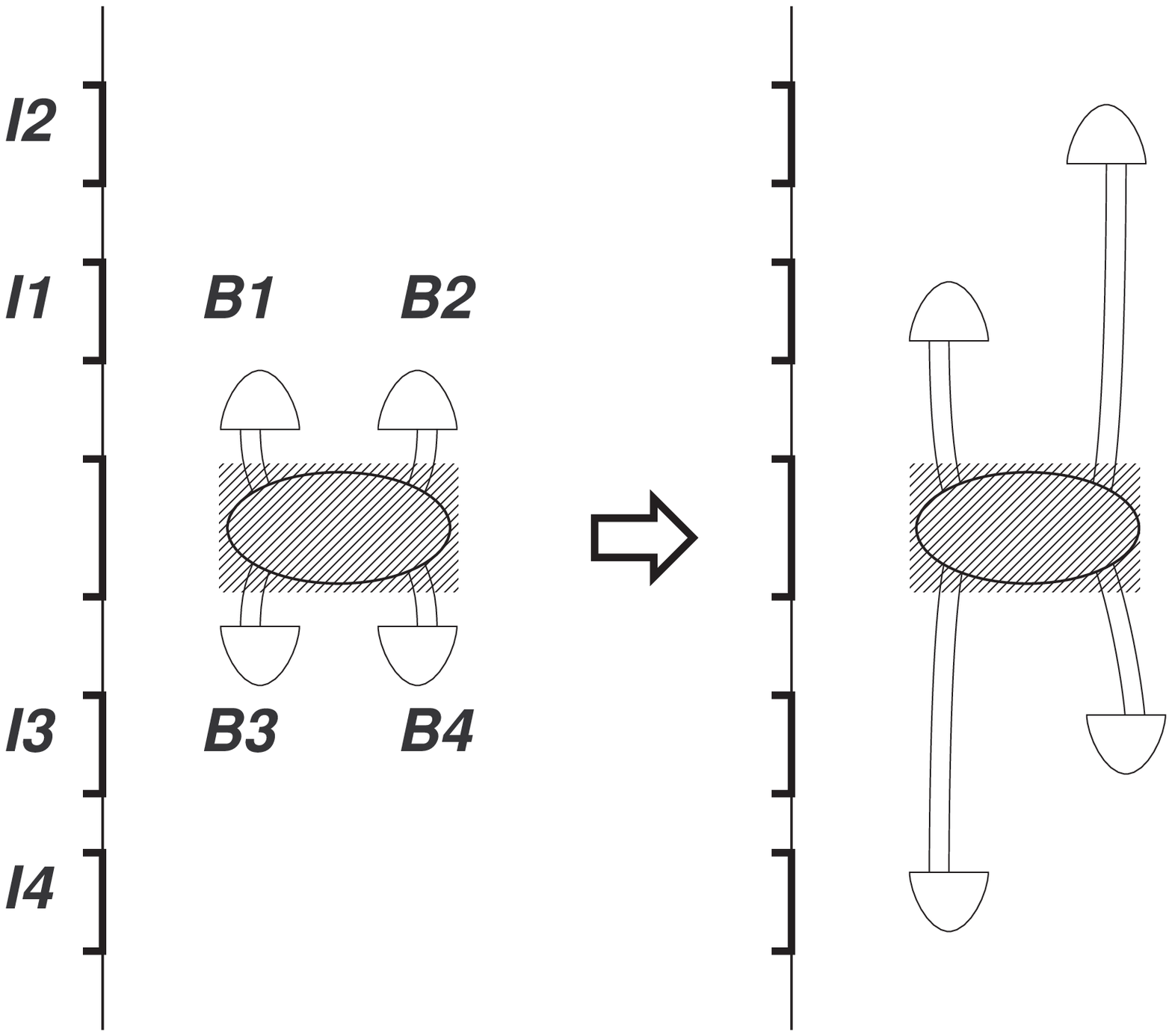}
\end{center}
\begin{center}
Figure 2.5
\end{center}
\end{figure}

\begin{enumerate}

\item 
$R_0$ is not altered. 

\item 
Each component of 
$L \setminus (R_0 \cup (\cup_{i=1}^p B_{s_i}))$ is preserved as a monotonic arc during the deformations.

\end{enumerate}

Then let $S_{t_1}, \dots , S_{t_q}$ be the bowl like 2-spheres which form the boundary of $C_0 \cup (\cup_{i=1}^p C_{s_i})$, and $B_{t_1}, \dots , B_{t_q}$ mutually disjoint 3-balls bounded by $S_{t_1}, \dots , S_{t_q}$. 
Then we apply the same kind of deformations as above, to make each $h( B_{t_i} )$ to be contained in $\text{Int} I_{t_i}$. 
By repeating such procedures finitely many times, we can obtain a desired isotopy $g_t$.\end{proof}

\section{Essential tangle decomposition from thin position}

In \cite{H-K}, we showed that for a link $L$ in a thin position, we can derive an essential tangle decomposition of $L$, which is closely related to the thin position. 
In this section, we quickly review the arguments, then describe some properties which will be required for the proof of the main result of this paper. 

\subsection{Deformation of tangle decomposition by bowl like 2-spheres}

We first describe a procedure for deforming a position of a link which admits a union of mutually disjoint bowl like 2-spheres giving a tangle decomposition with certain nice properties. 

In general, let $L$ be an unsplittable link, and 
$\mathcal{S} = S_1 \cup \dots \cup S_m$ a union of mutually disjoint bowl 
like 2-spheres. 
Let $C_0, C_1, \dots, C_m$ be the closures of the components of $S^3 \setminus \mathcal{S}$ such that $C_0$ lies exterior to all of $S_j$, and $C_i$ $(i=1, \dots, m)$ is the component lying directly inside of $S_i$. 
Then we suppose that these satisfy Property~1 in Section~2, and let $F_j$ $(j=0, 1, \dots , m)$ be as in Property~1. 

We say that a pair of arcs $\alpha_u$, $\alpha_b$ in $L \cap C_i$ is a {\it weakly bad pair of arcs} if the following conditions are satisfied. 

\begin{enumerate}
\item 
$\partial \alpha_u \subset F_i$ ($\partial \alpha_b \subset F_i$ resp.) and $\alpha_u$ ($\alpha_b$ resp.) contains exactly one critical point which is a maximum (minimum resp.). 

\item 
There exist arcs $\beta_u (\subset F_i)$, $\beta_b (\subset F_i)$, and disks $D_u (\subset S^3)$, $D_b (\subset S^3)$ such that $\partial D_u = \alpha_u \cup \beta_u$, $\partial D_b = \alpha_b \cup \beta_b$, $D_u\cap L=\alpha_u$, $D_b\cap L=\alpha_b$, $\text{Int} \beta_u \cap \text{Int} \beta_b = \emptyset$, $\text{Int} D_u \cap \text{Int} D_b = \emptyset$, $N( \beta_u, D_u )$ is above $F_i$, and $N( \beta_b, D_b )$ is below $F_i$.

\item
$\mid\alpha_u\cap\alpha_b\mid (=\mid\partial\alpha_u\cap\partial\alpha_b\mid) \ne 2$, i.e., that $\alpha_u\cup\alpha_b$ is not one component of $L$.
\end{enumerate}

\begin{remark}
We note that condition (3) does not, but should, appear in the definition of a weakly bad pair of arcs in \cite{H-K}.  Without this condition, Property~2 (following) can only guarantee that any weakly bad pair of arcs must form one component of $L$.  However, it is easy to verify that this (minor) mistake does not alter the other results of \cite{H-K}.
\end{remark}

Then we further suppose the following (see Property~2 in page 110 of \cite{H-K}).

\noindent 
{\bf 
Property~2}\qua
There does not exist a weakly bad pair of arcs for each $j=0, 1, \dots , m$. 

Suppose that there is a compressing disk $D$ for $\mathcal{S} \cap E(L)$ in $E(L)$. 
Let $S_k$ be the component of $\mathcal{S}$ such that $\partial D \subset S_k$ 
(hence, $C_k$ is directly inside of $S_k$), and let $C_l$ be the component which is directly outside of $S_k$. 
Since the argument is symmetric, we may suppose that $S_k$ is flat face up. 
Note that $D \subset C_l$ or $D \subset C_k$. 
Then we can show that Property~2 implies the following (see (*) in page 110 of \cite{H-K}).

\noindent 
{\bf 
Property~$\bf{2}'$}\qua
Suppose $D \subset C_l$ ($D \subset C_k$ resp.). 
Then all the critical points of $L$ in $C_l$ ($C_k$ resp.) are contained in a component of $C_l \setminus D$ ($C_k \setminus D$ resp.).

We call the closure of the component of $C_l \setminus D$ ($C_k \setminus D$ resp.) containing the level disk in the flat face of $S_k$ bounded by $\partial D$ {\it inside}, and call the closure of the other component {\it outside}. 

Then we have the following cases

\noindent
{\bf Case 1}\qua The disk $D$ is contained in $C_l$. 

We have the following subcases.

\noindent
{\bf Case 1.1}\qua No critical point of $L$ in $C_l$ is contained in the inside of $D$.

This case is divided into the following two cases.

\noindent
{\bf Case 1.1A}\qua 
No flat face up bowl like 2-sphere direcly inside of $S_l$ is contained in the inside of $D$.

In this case the inside of $D$ contains a collection of flat face down bowl like 2-spheres and monotonic arcs connecting $S_k$ and the flat face down bowl like 2-spheres. 
Then we pull down these interior bowl like 2-spheres into the interior of $C_k$, to obtain a new position of $L$ and a new system of bowl like 2-spheres. 
Note that this deformation of $L$ can be realized by an ambient isotopy of $S^3$, see \figref{3.1}.

\begin{figure}[ht]\small\anchor{3.1}
\begin{center}
\psfraga <-3pt,0pt> {D}{$D$}
\includegraphics[width=4.5cm, clip]{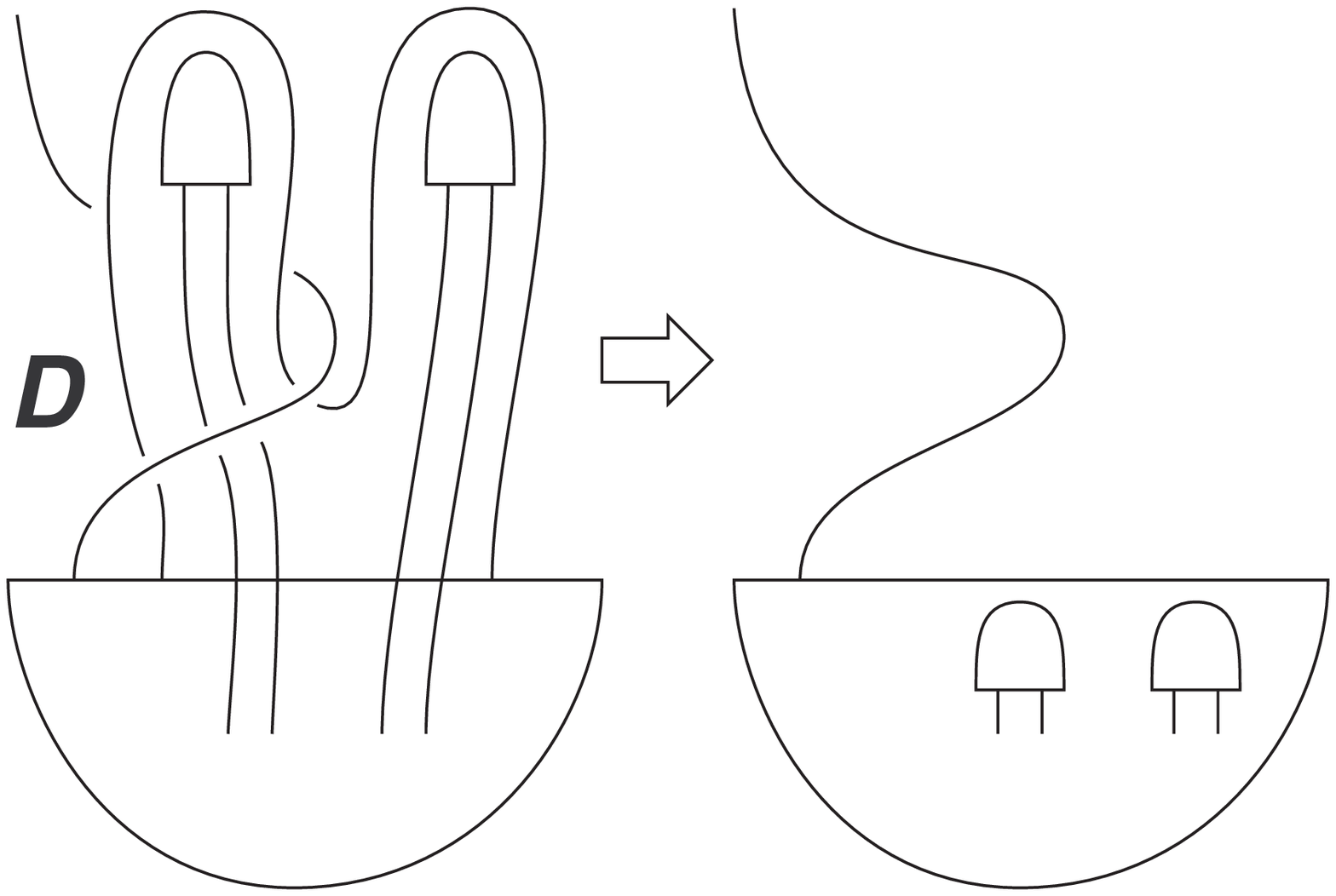}
\end{center}
\begin{center}
Figure 3.1
\end{center}
\end{figure}

\noindent
{\bf Case 1.1B}\qua 
There exists a flat face up bowl like 2-sphere component of $\partial C_l$ in the inside of $D$.

In this case, we first shrink the inside of $D$ similarly so that the inside of $D$ fits into a thin vertical cylinder above the flat face of $S_k$. 
Here we may suppose that the vertical cylinder is disjoint from $L \setminus N(\text{inside of }S_k)$. 
Now use a part of the cylinder together with a subdisk of a level 2-sphere and the subdisk of the flat face of $S_k$, to form a new flat face down bowl like 2-sphere. 
Then pull up the flat face down bowl like 2-sphere so that it is higher than the highest maximum of $L$ in $C_l$, and higher than flat face up bowl like 2-spheres in $(\text{outside of }D) \cap (\partial C_l \setminus S_l)$. 
Then we obtain a new position of $L$ and we take this flat face down bowl like 2-sphere together with the image of $\mathcal{S}$ under the above isotopy as a new system of bowl like 2-spheres, see \figref{3.2}.

\begin{figure}[ht]\small\anchor{3.2}
\begin{center}
\includegraphics[width=6cm, clip]{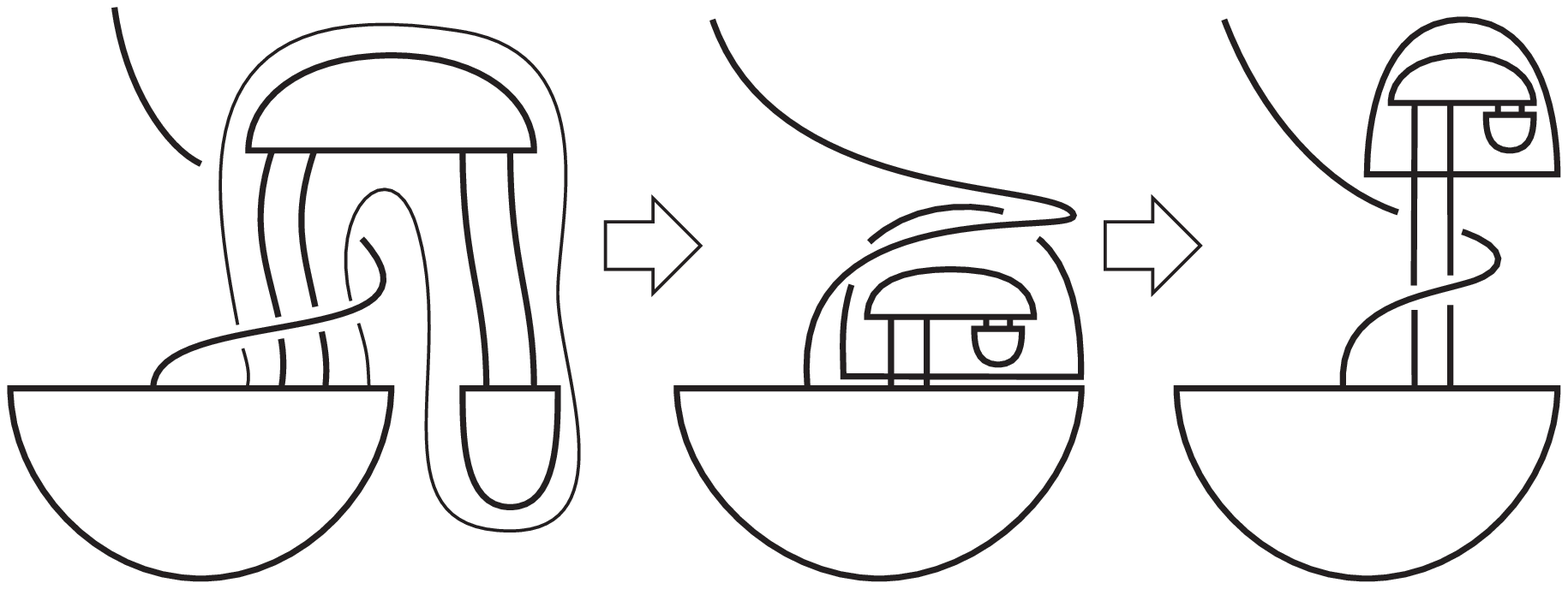}
\end{center}
\begin{center}
Figure 3.2
\end{center}
\end{figure}

\noindent
{\bf Case 1.2}\qua
All critical points of $L$ in $C_l$ are contained in the inside of $D$.

This case is divided into the following two cases.

\noindent
{\bf Case 1.2A}\qua 
No flat face up bowl like 2-sphere in $\partial C_l$ is in the inside of $D$.

\noindent
{\bf Case 1.2B}\qua 
There exists a flat face up bowl like 2-sphere in $\partial C_l$ in the inside of $D$.

In both cases, we apply the argument of Case~1.1B to obtain a new position of $L$ and a system of bowl like 2-spheres. 
See Figures~\fref{3.3} and \fref{3.4}.

\begin{figure}[ht]\small
 \begin{minipage}{0.5\hsize}\anchor{3.3}
  \begin{center}
   \includegraphics[width=5cm, clip]{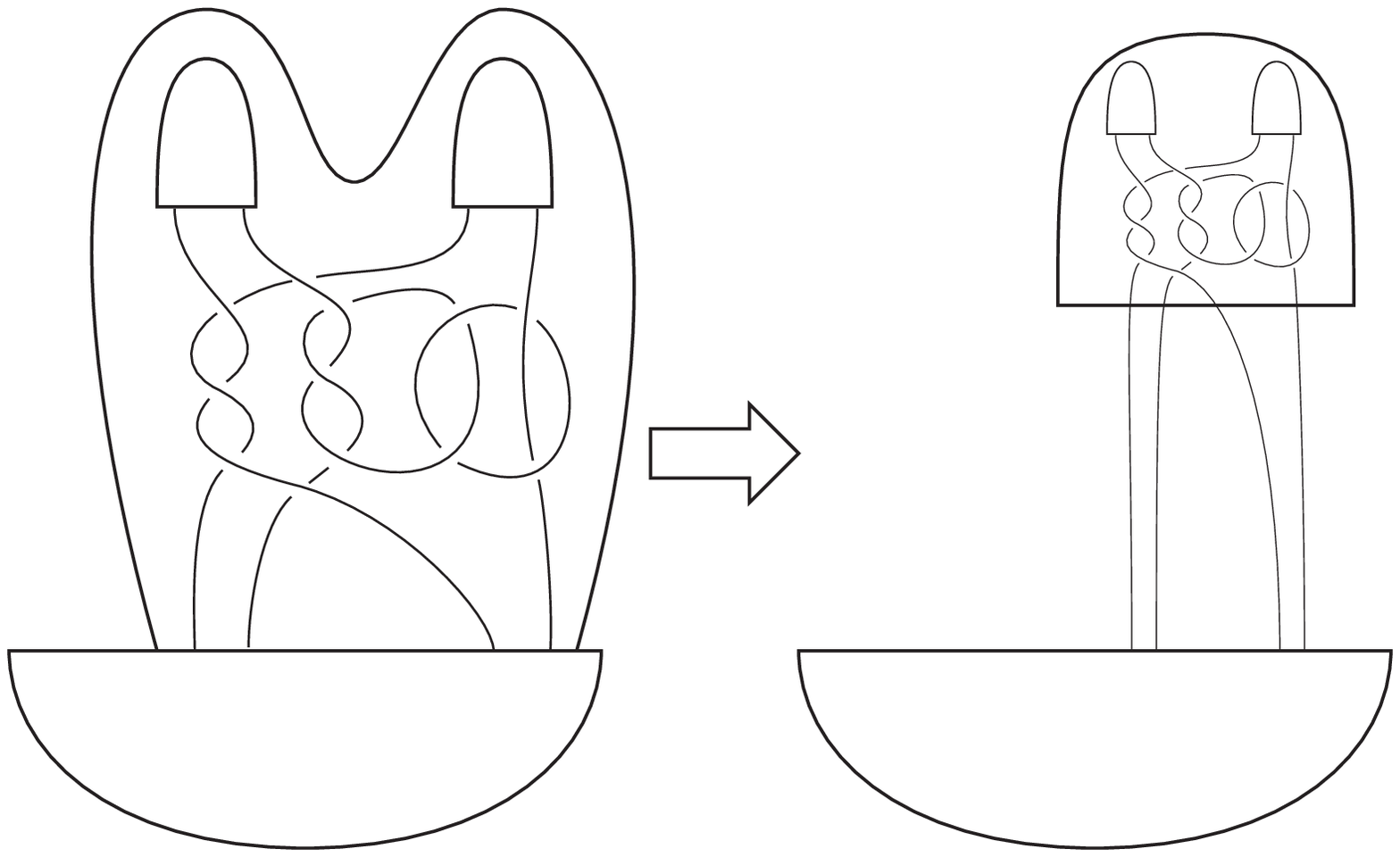}
   
   Figure 3.3
  \end{center}

 \end{minipage}
 \begin{minipage}{0.5\hsize}\anchor{3.4}
  \begin{center}
   \includegraphics[width=5cm, clip]{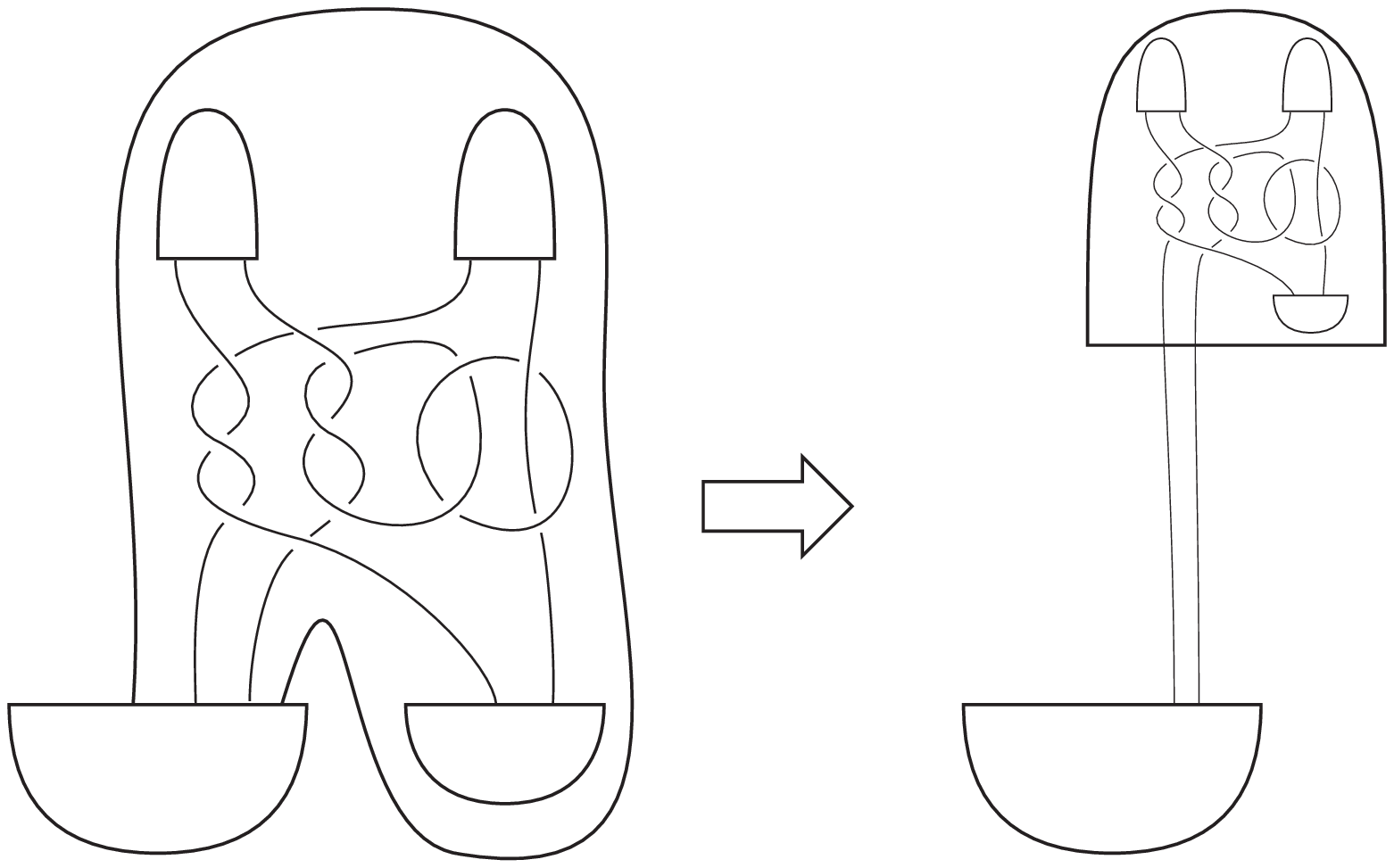}
   
   Figure 3.4
  \end{center}
 \end{minipage}
\end{figure}

\noindent 
{\bf Case 2}\qua 
The disk $D$ is contained in $C_k$.

This case is divided into the following four subcases.

\noindent
{\bf Case 2.1A}\qua
No critical point of $L$ in $C_k$ is contained in the inside of $D$, and no flat face down bowl like 2-sphere in $\partial C_k$ is in the inside of $D$, see \figref{3.5.1}.

\noindent
{\bf Case 2.1B}\qua 
No critical point of $L$ in $C_k$ is contained in the inside of $D$, and there exists a flat face down bowl like 2-sphere in $\partial C_k$ in the inside of $D$, see \figref{3.5.2}.

\noindent
{\bf Case 2.2A}\qua 
All critical points of $L$ in $C_k$ are contained in the inside of $D$, and no flat face down bowl like 2-sphere in $\partial C_k$ is in the inside of $D$, see \figref{3.5.3}.

\noindent
{\bf Case 2.2B}\qua 
All critical points of $L$ in $C_k$ are contained in the inside of $D$, and there exists a flat face down bowl like 2-sphere in $\partial C_k$ in the inside of $D$, see \figref{3.5.4}.

Each of these cases are handled analogously to the corresponding subcases of Case~1.

\begin{figure}[ht]\small
 \begin{minipage}{0.5\hsize}\anchor{3.5.1}
  \begin{center}
   \includegraphics[width=5cm, clip]{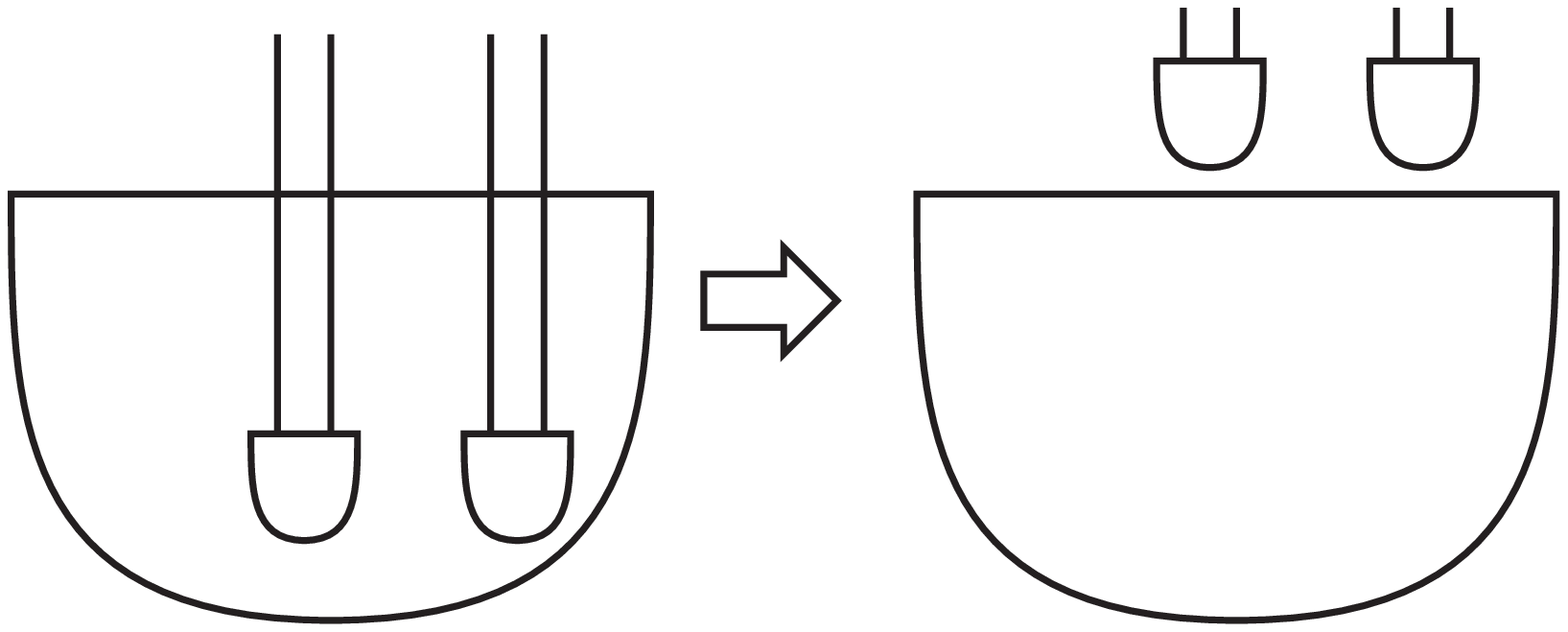}
   
   Figure 3.5.1
  \end{center}

 \end{minipage}
 \begin{minipage}{0.5\hsize}\anchor{3.5.2}
  \begin{center}
   \includegraphics[width=5cm, clip]{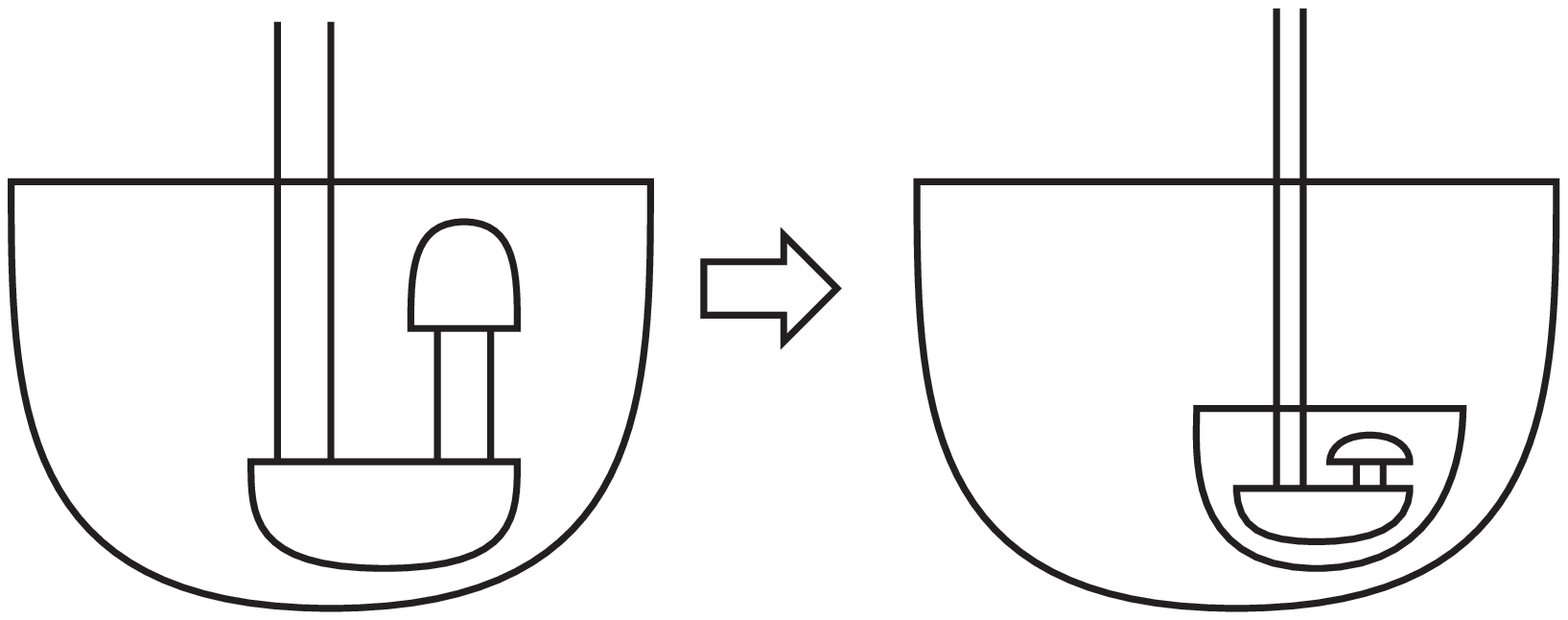}
   
   Figure 3.5.2
  \end{center}
 \end{minipage}
\end{figure}

\begin{figure}[ht]\small
 \begin{minipage}{0.5\hsize}\anchor{3.5.3}
  \begin{center}
   \includegraphics[width=5cm, clip]{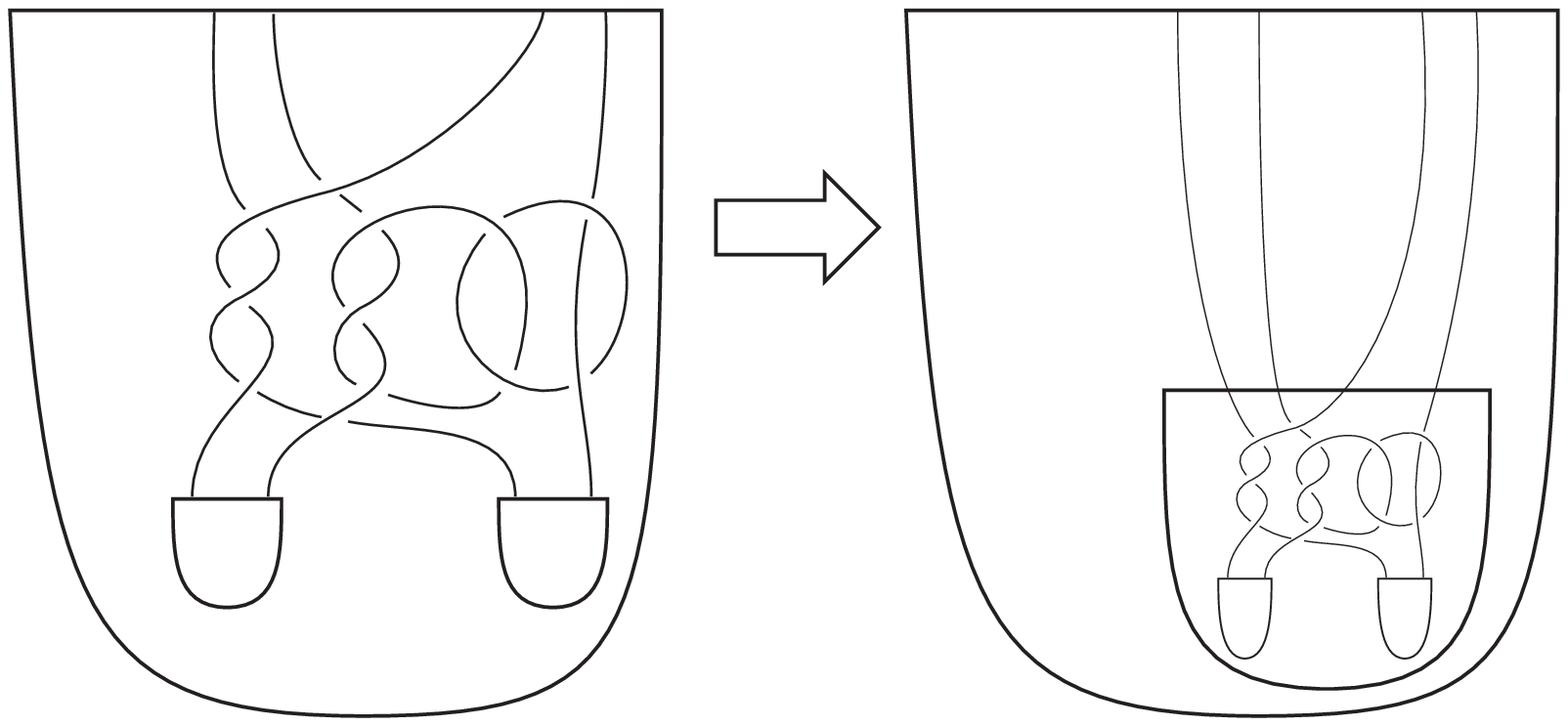}
   
   Figure 3.5.3
  \end{center}

 \end{minipage}
 \begin{minipage}{0.5\hsize}\anchor{3.5.4}
  \begin{center}
   \includegraphics[width=5cm, clip]{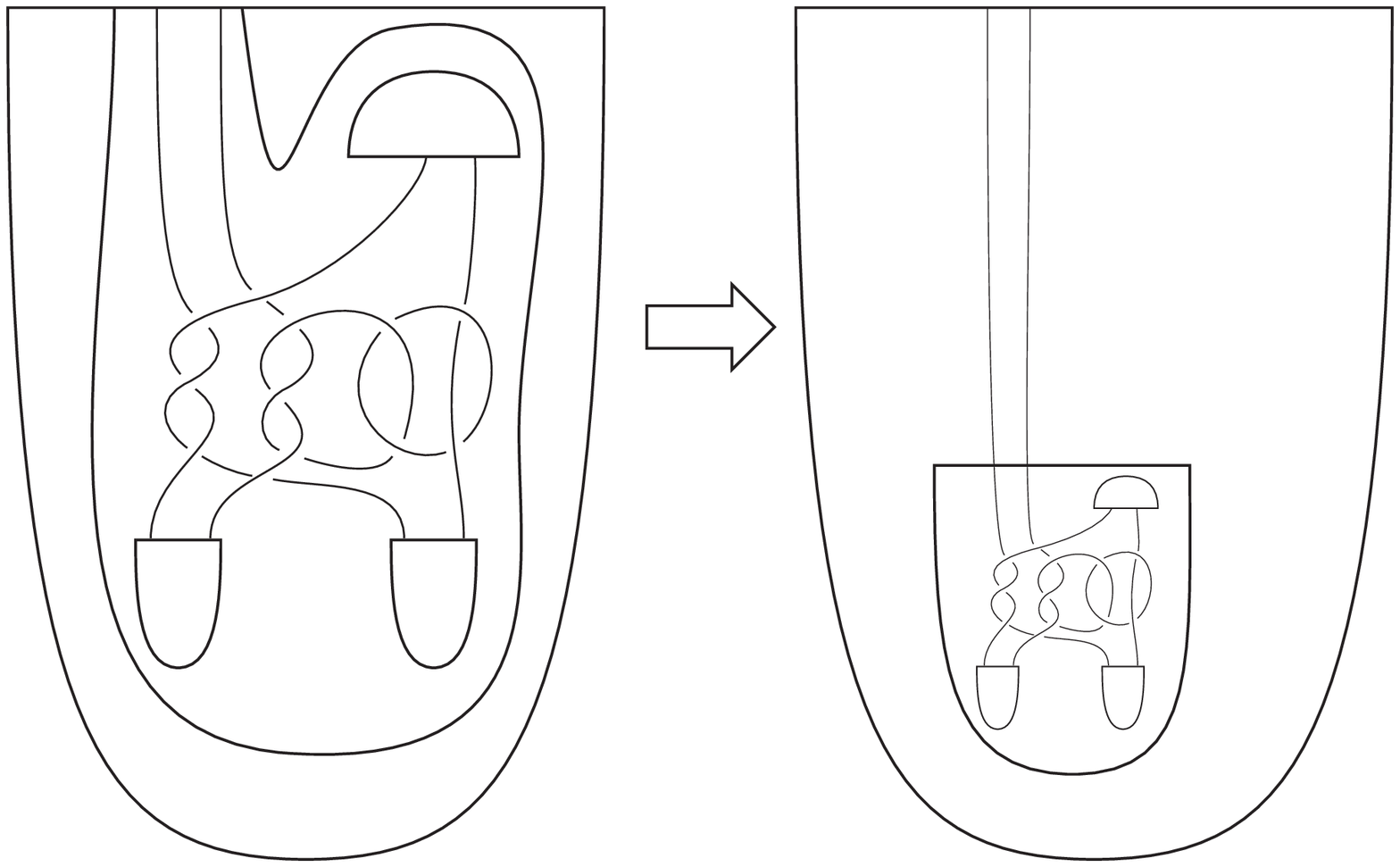}
   
   Figure 3.5.4
  \end{center}
 \end{minipage}
\end{figure}

\subsection{From a thin position to an essential tangle decomposition}

Now we describe the procedures in \cite{H-K}. 
Let $L (= L^{(0)})$ be a link in a thin position together with a union of bowl like 2-spheres $S^{(0)} = S_1^{(0)} \cup \dots \cup S_{m_0}^{(0)}$ obtained as in the last paragraph in page 108 of \cite{H-K}.
Note that $S^{(0)}$ here refers to a collection of bowl like 2-spheres coming from thin position of $L$, while $\mathcal{S}$ in previous sections refered to a general collection of bowl like 2-spheres not directly related to thin position.
Suppose there exists a compressing disc $D$ for $S^{(0)}\cap E(L)$ in $E(L)$. 
We showed that $S^{(0)}$ satisfies Property~1 of Section~2 and Property~2 (hence Property~$2'$) of Section~3.1. 
Hence, we can apply one of the above procedures to $L^{(0)}$, $S^{(0)}$ to obtain a new position of the link, say $L^{(1)}$, and a union of bowl like 2-spheres $S^{(1)} = S_1^{(1)} \cup \dots \cup S_{m_1}^{(1)}$ giving a tangle decomposition of $L^{(1)}$. 
If $S^{(1)} \cap E( L^{(1)} )$ is compressible in $E( L^{(1)} )$, then we can show that $L^{(1)}$, $S^{(1)}$ also satisfies Properties~1,2. 
Hence we can apply one of the above procedures to $L^{(1)}$, $S^{(1)}$ to obtain $L^{(2)}$, and $S^{(2)} = S_1^{(2)} \cup \dots \cup S_{m_2}^{(2)}$. 
By using what we call Scharlemann-Thompson complexity for $L^{(i)}$, $S^{(i)}$, we showed that the procedure terminates in finitely many steps to give a position $L^{(n)}$ of $L$, and a union of mutually disjoint bowl like 2-spheres $S^{(n)} = S_1^{(n)} \cup \dots \cup S_{m_n}^{(n)}$ giving an essential tangle decomposition, that is to say $S^{(n)}\cap E(L)$ are essential in $E(L)$. 
Note that $L^{(n)}$, $S^{(n)}$ satisfies Property~1. 

\subsection{Cocoons}

Let $L^{(k)}$, $S^{(k)} = S_1^{(k)} \cup \dots \cup S_{m_k}^{(k)}$ $(k=1, \dots , n)$ be as in Section~3.2. 
Obviously, we can retrieve the original thin position $L (= L^{(0)})$ by tracing the sequence of the ambient isotopies conversely. 
In this subsection, we will make an observation about how the cocoons of $L^{(n)}$, $S^{(n)}$ survive in $L^{(0)}$, $S^{(0)}$. 
Let $C_0^{(k)}$, $C_1^{(k)}, \dots , C_{m_k}^{(k)}$ be the closures of the components of $S^3 \setminus S^{(k)}$ such that $C_0^{(k)}$ lies exterior to all of the $S_i^{(k)}$, and $C_i^{(k)}$ $(i=1, \dots , m_k)$ is the component lying directly inside of $S_i^{(k)}$. 
Let $R_j^{(k)}$ be the cocoon of $L^{(k)}$ in $C_j^{(k)}$. 
Recall that $h: S^2 \times {\mathbb R} \rightarrow {\mathbb R}$ is the height function. 
We first show the following.

\noindent
{\bf 
Assertion 1 
}
We may suppose that $h( R_0^{(0)} ), h( R_1^{(0)} ), \dots, h( R_{m_0}^{(0)} )$ are mutually disjoint intervals in ${\mathbb R}$. 

\begin{proof}
If there is a critical point of $L$ in $C_i^{(0)}$, then we let $\text{Min}_i^{(0)}$ ($\text{Max}_i^{(0)}$ resp.) be the height of the lowest minimum (the highest maximum resp.) of $L^{(0)}$ in $R_i^{(0)}$ (see (1) of Remark~1 in Section~2). 
If there is no critical point of $L^{(0)}$ in $R_i^{(0)}$, then we let $\text{Min}_i^{(0)}$ ($\text{Max}_i^{(0)}$ resp.) be the minimum (maximum resp.) of $h( R_i^{(0)} )$. 
For the proof of Assertion~1, it is enough to show that for each pair $i, j$ ($i \ne j$), we have $[ \text{Min}_i^{(0)}, \text{Max}_i^{(0)} ] \cap [ \text{Min}_j^{(0)}, \text{Max}_j^{(0)} ] = \emptyset$. 
Suppose that $[ \text{Min}_i^{(0)}, \text{Max}_i^{(0)} ] \cap [ \text{Min}_j^{(0)}, \text{Max}_j^{(0)} ] \ne \emptyset$ for some pair $i, j$. 
We divide the proof into the following cases.

\noindent
{\bf Case 1}\qua 
For any pair $i, j$ ($i \ne j$) with 
$[ \text{Min}_i^{(0)}, \text{Max}_i^{(0)} ] \cap [ \text{Min}_j^{(0)}, \text{Max}_j^{(0)} ] \ne \emptyset$, there does not exist a critical point of $L^{(0)}$ in $R_i^{(0)} \cup R_j^{(0)}$.

In this case, it is easy to see that we can shrink the cocoons under consideration vertically without changing the position of $L^{(0)}$. 
By applying slight isotopy if necessary, we can make the images of them by $h$ to be mutually disjoint, see \figref{3.6}.

\begin{figure}[ht]\small\anchor{3.6}
\begin{center}
\includegraphics[width=8cm, clip]{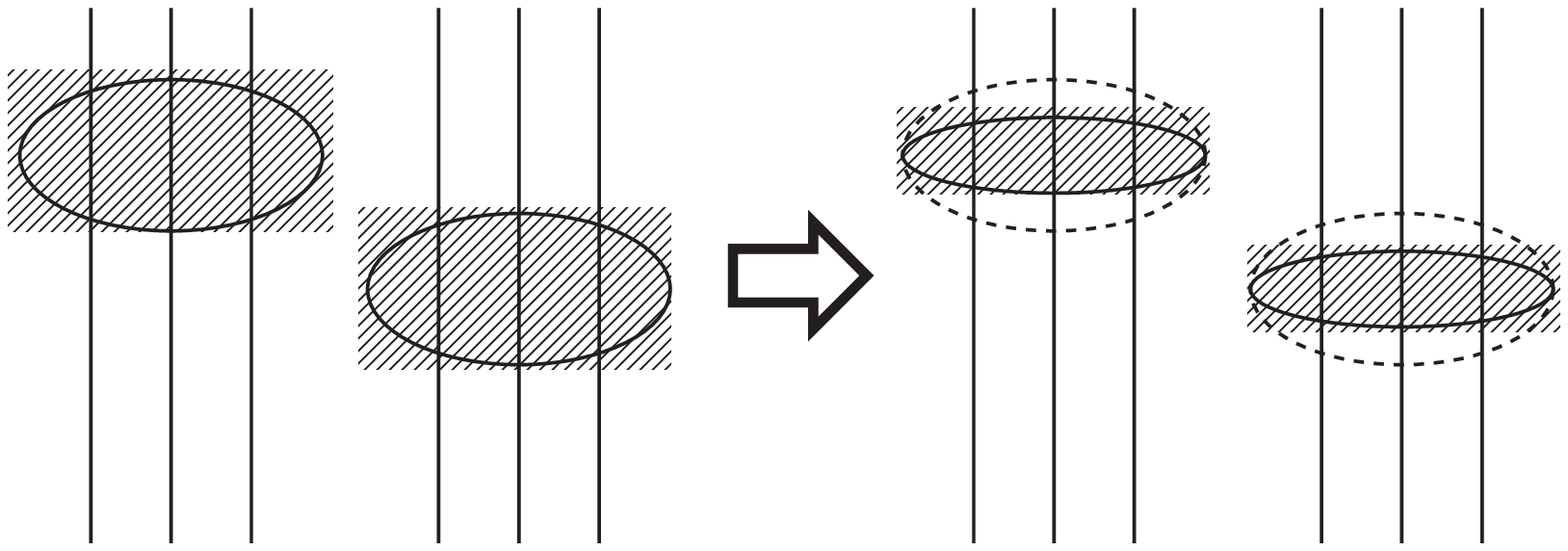}
\end{center}
\begin{center}
Figure 3.6
\end{center}
\end{figure}

\noindent
{\bf Case 2}\qua
There exists a pair $i, j$ ($i \ne j$) with $[ \text{Min}_i^{(0)}, \text{Max}_i^{(0)} ] \cap [ \text{Min}_j^{(0)}, \text{Max}_j^{(0)} ]$ $\ne \emptyset$ such that $R_i^{(0)}$ contains a minimum (hence maximum too) of $L^{(0)}$. 
This case is divided into the following two subcases.

\noindent
{\bf Case 2.1}\qua 
There does not exist a pair $i, j$ ($i \ne j$) with $[ \text{Min}_i^{(0)}, \text{Max}_i^{(0)} ] \cap [ \text{Min}_j^{(0)}, \text{Max}_j^{(0)} ] \ne \emptyset$ such that both $R_i^{(0)}$, $R_j^{(0)}$ contain critical points of $L^{(0)}$.

In this case, the argument for Case~1 basically works. 
The only difference is that if we apply vertical shrink to the cocoons containing critical points, then the position of $L^{(0)}$ is changed (however, it is clear that the width is unchanged), see \figref{3.7}.

\begin{figure}[ht]\small\anchor{3.7}
\begin{center}
\includegraphics[width=8cm, clip]{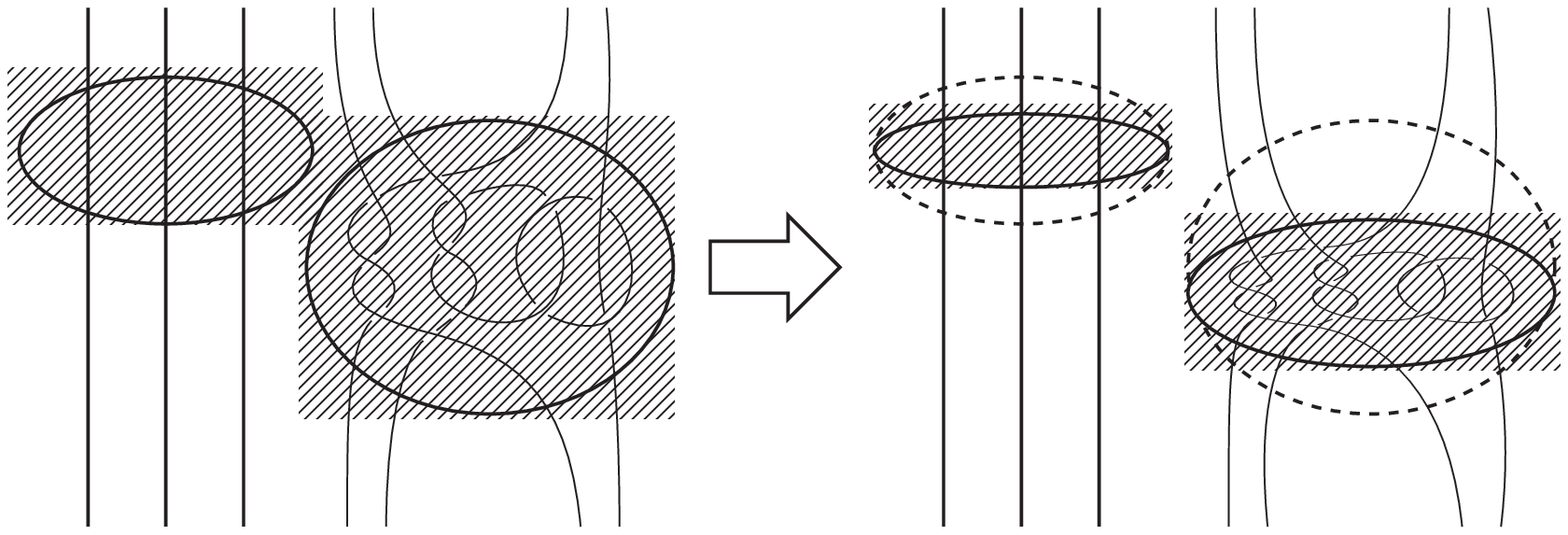}
\end{center}
\begin{center}
Figure 3.7
\end{center}
\end{figure}

\noindent
{\bf Case 2.2}\qua 
There exists a pair $i, j$ ($i \ne j$) with $[ \text{Min}_i^{(0)}, \text{Max}_i^{(0)} ] \cap [ \text{Min}_j^{(0)}, \text{Max}_j^{(0)} ]$ $\ne \emptyset$ such that both $R_i^{(0)}$, $R_j^{(0)}$ contain critical points of $L^{(0)}$.

Without loss of generality, we may suppose that $\text{Max}_j^{(0)}$ or $\text{Min}_j^{(0)}$ is contained in $[ \text{Min}_i^{(0)}, \text{Max}_i^{(0)} ]$. 
Let $\theta \in [ \text{Min}_i^{(0)}, \text{Max}_i^{(0)} ]$ be a height such that $\vert h^{-1}( \theta ) \cap L^{(0)} \vert$ is minimal among all points in $[ \text{Min}_i^{(0)}, \text{Max}_i^{(0)} ]$. 
Then we shrink $R_i^{(0)}$ vertically so that $h( R_i^{(0)})$ is a very small neighborhood of $\theta$. 
This new position of $L^{(0)}$ is obviously thinner than that of $L^{(0)}$, a contradiction, see \figref{3.8}.

\begin{figure}[ht]\small\anchor{3.8}
\begin{center}
\includegraphics[width=8cm, clip]{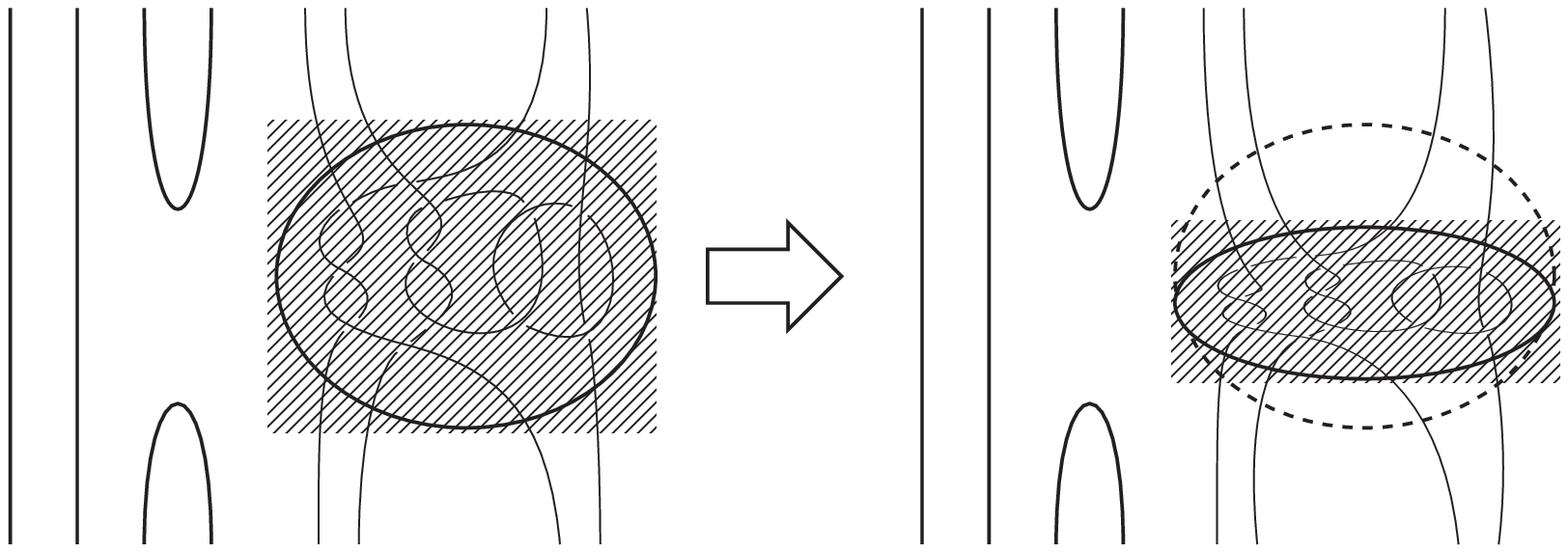}
\end{center}
\begin{center}
Figure 3.8
\end{center}
\end{figure}

This completes the proof of Assertion~1. 
\end{proof}

\noindent
{\bf How cocoons $\{ R_j^{(k)}\}$ are related to $\{ R_j^{(k+1)}\}$} 

We next make an observation about how the cocoons $\{ R_j^{(k)}\}$ are related to $\{ R_j^{(k+1)}\}$. 
We analyze the situation by cases as in 3.1.

\noindent
{\bf Case 1.1A}\qua 
Suppose that $L^{(k+1)}$, $S^{(k+1)}$ are obtained from $L^{(k)}$, $S^{(k)}$ by applying the deformation of Case~1.1A. 
We can retrieve $L^{(k)}$, $S^{(k)}$ from  $L^{(k+1)}$, $S^{(k+1)}$ as follows. 

In this case, $m_{k+1} = m_k$. 
We may suppose that $S_j^{(k)}$ $(j=1, \dots , m_k)$ corresponds to $S_j^{(k+1)}$. 
Suppose that the deformation from $L^{(k)}$, $S^{(k)}$  to $L^{(k+1)}$, $S^{(k+1)}$ is performed by using a compression disk for $S_u^{(k)}$ (recall that $C_u^{(k)}$ is directly inside of $S_u^{(k)}$), and let $C_v^{(k)}$ be the component directly outside of $S_u^{(k)}$. 
Let $\widetilde{S^{(k)}}$ be the union of the flat face down bowl like 2-spheres that are pulled into $C_u^{(k)}$, and $\widetilde{S^{(k+1)}}$ the corresponding 2-spheres in $S^{(k+1)}$. 
Here $\widetilde{S^{(k+1)}}$ is obtained from $\widetilde{S^{(k)}}$ by applying similarity deformations and then parallel translations for each component. 

Then the  position $L^{(k)}$ is obtained from $L^{(k+1)}$ as follows. 
First let $\widetilde{B^{(k+1)}}$ be the union of mutually disjoint 3-balls bounded by $\widetilde{S^{(k+1)}}$. 
We take the piece $L^{(k+1)} \cap \widetilde{B^{(k+1)}}$, and deform it by a similarity deformation and parallel translation (which is in fact the inverse of the above deformation) to put the piece in a position higher than $R_v^{(k+1)}$. 
Note that in this stage some components of $L^{(k+1)} \cap \widetilde{B^{(k+1)}}$ are torn into 
a union of arcs. 
Then we add monotonic arcs to obtain $L^{(k)}$. 

This observation implies the following.

\noindent
{\bf 
Facts 1.1A
}
\begin{enumerate}
\item 
If $i \ne v$, then the pairs $(R_i^{(k)}, L^{(k)} \cap R_i^{(k)})$ and $(R_i^{(k+1)}, L^{(k+1)} \cap R_i^{(k+1)})$ are similar. 

\item 
$(R_v^{(k)}, L^{(k)} \cap R_v^{(k)})$ is obtained from $(R_v^{(k+1)}, L^{(k+1)} \cap R_v^{(k+1)})$ by adding some monotonic arcs. 
\end{enumerate}

See \figref{3.9}. 

\begin{figure}[ht]\small\anchor{3.9}
\begin{center}
\psfraga <0pt,1pt> {R1}{$R_v^{(k+1)}$}
\psfraga <0pt,-1pt> {R2}{$R_u^{(k+1)}$}
\includegraphics[width=9cm, clip]{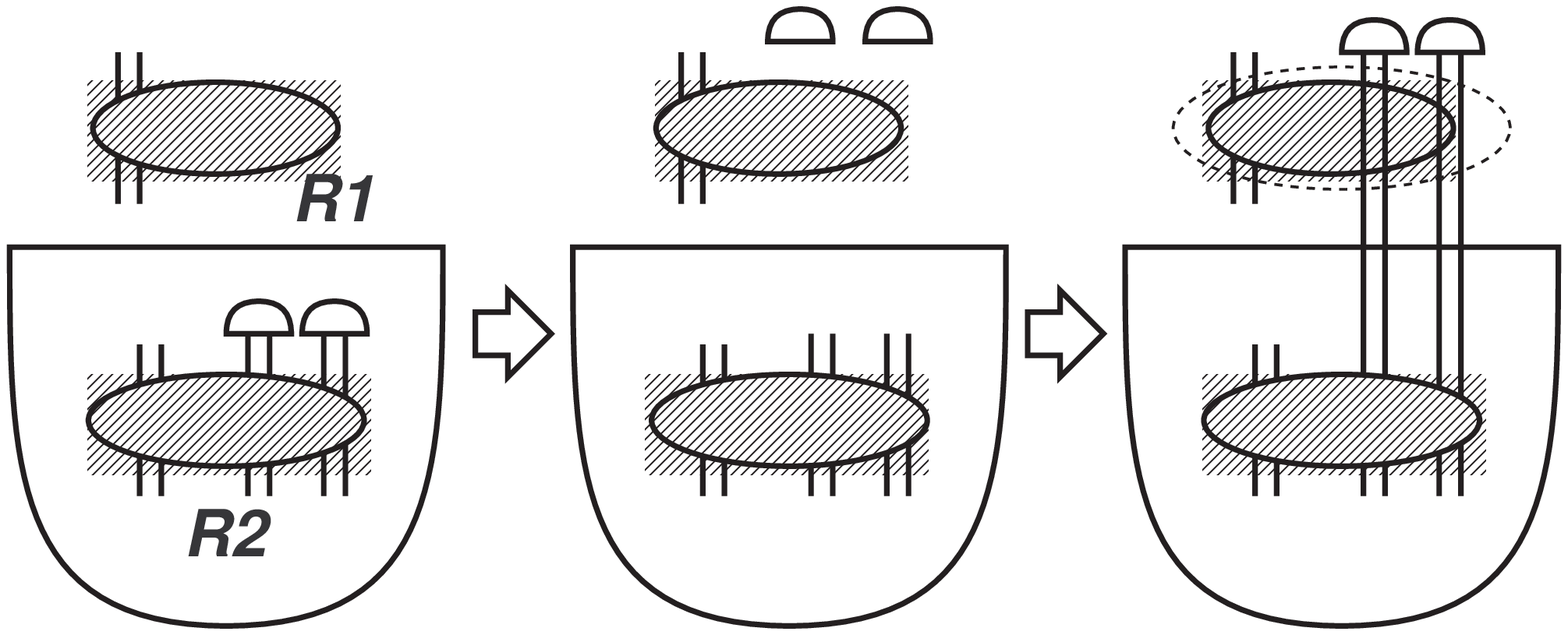}
\end{center}
\begin{center}
Figure 3.9
\end{center}
\end{figure}

\noindent
{\bf Case 1.1B}\qua 
In this case $m_{k+1} = m_k +1$. 
We may suppose that $S_j^{(k)}$ $(j=1, \dots , m_k)$ corresponds to $S_j^{(k+1)}$, and $S_{m_{k+1}}^{(k+1)}$ is the extra component. 
Let $C_u^{(k+1)}$ be the component which is directly outside of $S_{m_{k+1}}^{(k+1)}$. 
Let $\widetilde{S_d^{(k+1)}}$ ($\widetilde{S_u^{(k+1)}}$ resp.) be the union of the components of $\partial C_{m_{k+1}}^{(k+1)} \setminus S_{m_{k+1}}^{(k+1)}$ which are flat face down (up resp.) bowl like 2-spheres. 
Then $L^{(k)}$ is obtained from $L^{(k+1)}$ as follows. 
Let $\widetilde{B_d^{(k+1)}}$ ($\widetilde{B_u^{(k+1)}}$ resp.) be the union of mutually disjoint 3-ball(s) bounded by $\widetilde{S_d^{(k+1)}}$ ($\widetilde{S_u^{(k+1)}}$ resp.). 
We take the piece $L^{(k+1)} \cap \widetilde{B_d^{(k+1)}}$ ($L^{(k+1)} \cap \widetilde{B_u^{(k+1)}}$ resp.) and deform it by a similarity deformation and parallel translation to put the piece in a position higher (lower resp.) than $R_u^{(k+1)}$. 
Note that in this stage some components of 
$L^{(k+1)} \cap (\widetilde{B_d^{(k+1)}} \cup \widetilde{B_u^{(k+1)}})$ 
are torn into a union of arcs. 
Then we replace some monotonic arcs in $L^{(k+1)} \cap C_u^{(k+1)}$
with other monotonic arcs to obtain $L^{(k)}$.

This implies the following.

\noindent
{\bf 
Facts 1.1B
}
\begin{enumerate}
\item 
If $i \ne u, m_{k+1}$, then the pairs $(R_i^{(k)}, L^{(k)} \cap R_i^{(k)})$ and $(R_i^{(k+1)}, L^{(k+1)} \cap R_i^{(k+1)})$ are similar. 

\item 
Note that $L^{(k+1)} \cap R_{m_{k+1}}^{(k+1)}$ consists of monotonic arcs. 
This fact allows us to regard that $(R_u^{(k)}, L^{(k)} \cap R_u^{(k)})$ is obtained from $R_u^{(k+1)}$ and $R_{m_{k+1}}^{(k+1)}$ by putting them in a vertically disjoint position ($R_u^{(k+1)}$ is below, and $R_{m_{k+1}}^{(k+1)}$ is above) and adding some monotonic arcs. 
\end{enumerate}

See \figref{3.10}.

\begin{figure}[ht]\small\anchor{3.10}
\begin{center}
\psfraga <0pt,-1pt> {R1}{$R_{m_{k+1}}^{(k+1)}$}
\psfraga <0pt,-1pt> {R2}{$R_u^{(k+1)}$}
\includegraphics[width=8cm, clip]{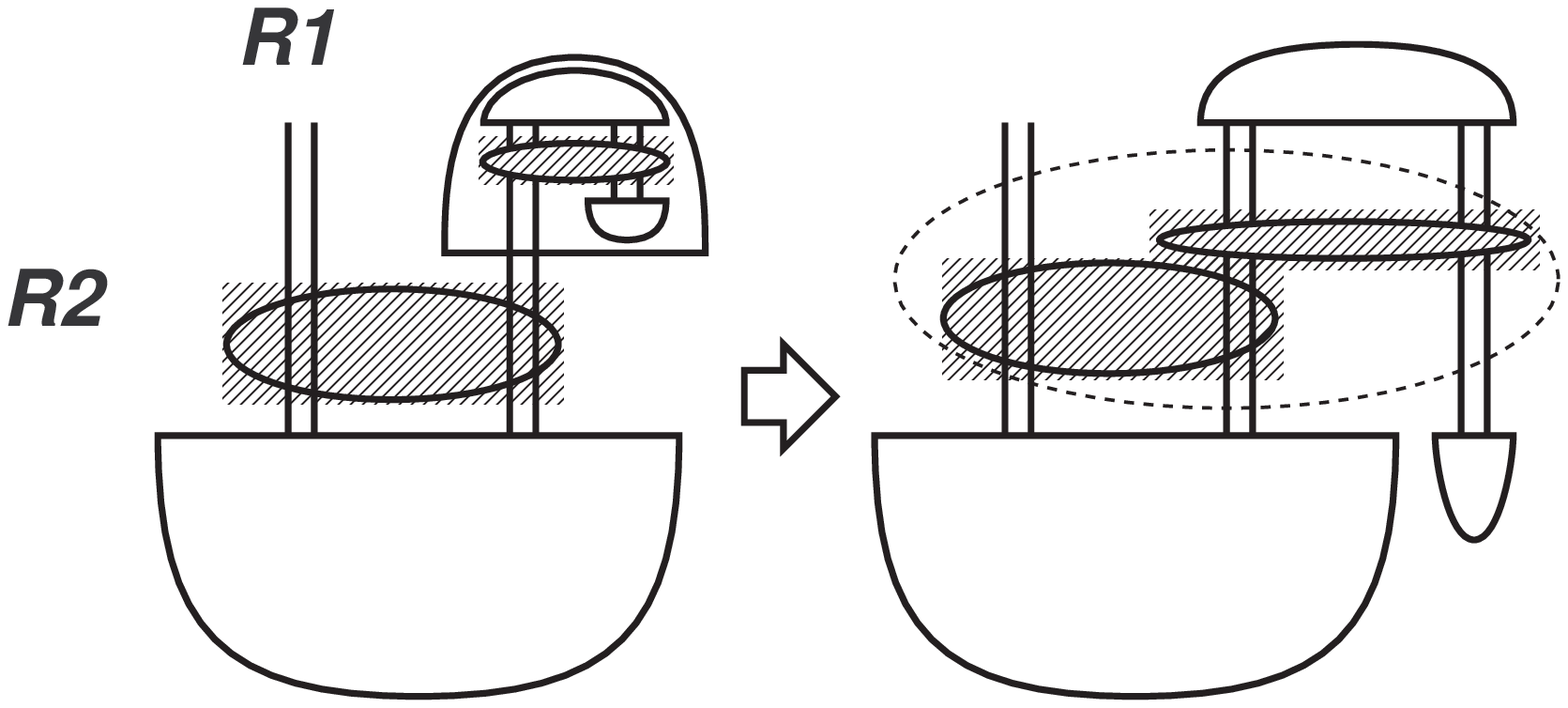}
\end{center}
\begin{center}
Figure 3.10
\end{center}
\end{figure}

\noindent
{\bf Case 1.2A}\qua 
In this case $m_{k+1} = m_k +1$. 
We may suppose that $S_j^{(k)}$ $(j=1, \dots , m_k)$ corresponds to $S_j^{(k+1)}$, and $S_{m_{k+1}}^{(k+1)}$ is the extra component. 
Let $C_u^{(k+1)}$ be the component which is directly outside of $S_{m_{k+1}}^{(k+1)}$. 
Let $\widetilde{S^{(k+1)}} = \partial C_{m_{k+1}}^{(k+1)} \setminus S_{m_{k+1}}^{(k+1)}$.
Then $L^{(k)}$ is obtained from $L^{(k+1)}$ as follows. 
Let $\widetilde{B^{(k+1)}}$ be the union of mutually disjoint 3-balls bounded by $\widetilde{S^{(k+1)}}$. 
We first remove the piece $L^{(k+1)} \cap (C_u^{(k+1)} \cup C_{m_{k+1}}^{(k+1)} \cup \widetilde{B^{(k+1)}})$ from $L^{(k+1)}$. 
Then add the pieces $L^{(k+1)} \cap \widetilde{B^{(k+1)}}$ and $L^{(k+1)} \cap R_{m_{k+1}}^{(k+1)}$ by applying similarity deformations and parallel translations into a position appropriate for retrieving $L^{(k)}$. 
Finally add some monotonic arcs to obtain $L^{(k)}$. 

This observation implies the following.

\noindent
{\bf 
Facts 1.2A 
}
\begin{enumerate}
\item 
If $i \ne u, m_{k+1}$, then the pairs $(R_i^{(k)}, L^{(k)} \cap R_i^{(k)})$ and $(R_i^{(k+1)}, L^{(k+1)} \cap R_i^{(k+1)})$ are similar. 

\item 
Note that $L^{(k+1)} \cap R_u^{(k+1)}$ consists of monotonic arcs. 
This fact allows us to regard $(R_u^{(k)}, L^{(k)} \cap R_u^{(k)})$ as obtained from $R_u^{(k+1)}$ and $R_{m_{k+1}}^{(k+1)}$ by putting them in a vertically disjoint position ($R_u^{(k+1)}$ is below, and $R_{m_{k+1}}^{(k+1)}$ is above) and adding some monotonic arcs. 
\end{enumerate}

See \figref{3.11}.

\begin{figure}[ht]\small\anchor{3.11}
\begin{center}
\psfraga <-2pt,-1pt> {R1}{$R_{m_{k+1}}^{(k+1)}$}
\psfraga <0pt,-1pt> {R2}{$R_u^{(k+1)}$}
\includegraphics[width=8cm, clip]{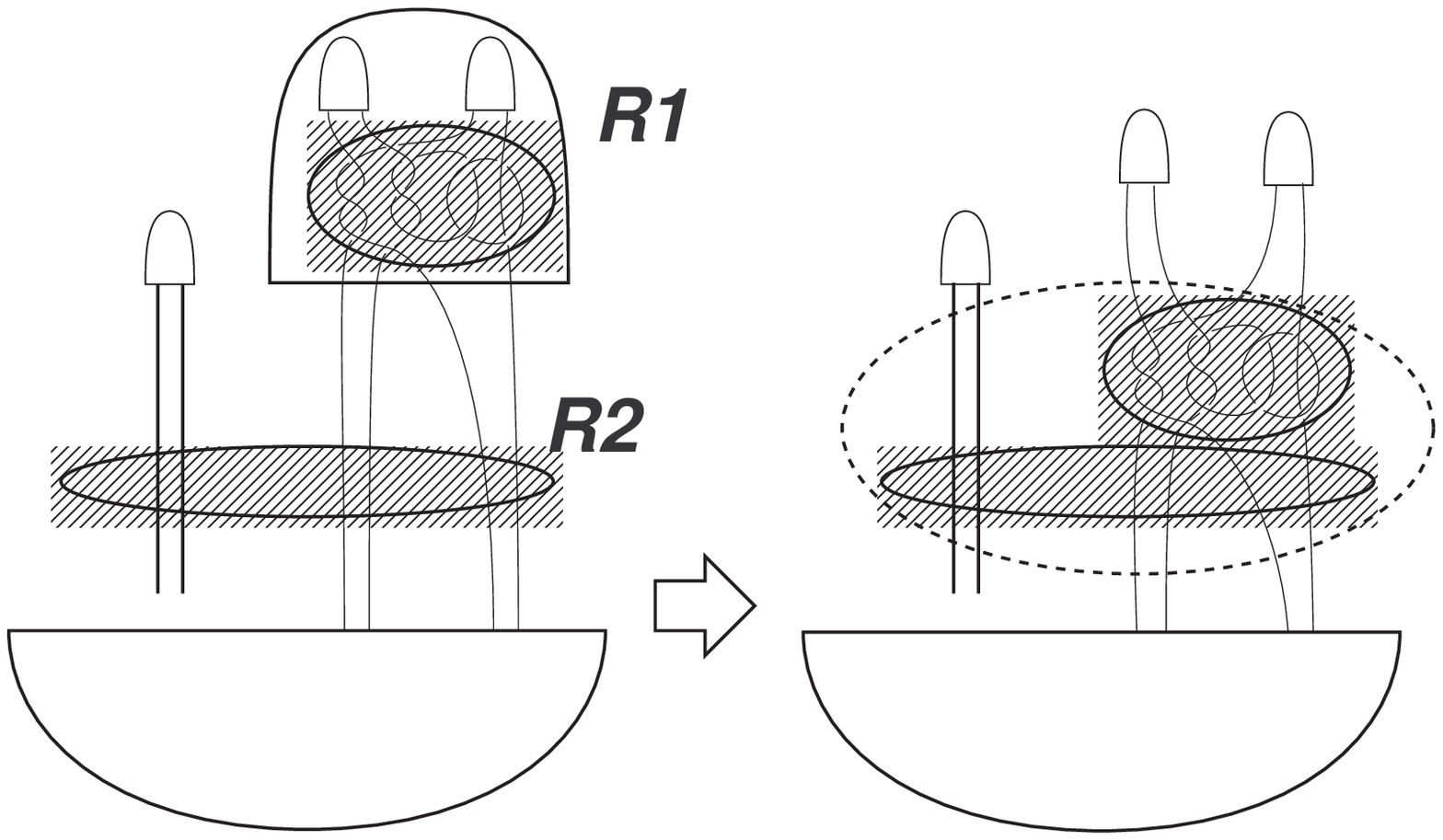}
\end{center}
\begin{center}
Figure 3.11
\end{center}
\end{figure}

\noindent
{\bf Case 1.2B}\qua 
In this case $m_{k+1} = m_k +1$. 
We may suppose that $S_j^{(k)}$ $(j=1, \dots , m_k)$ corresponds to $S_j^{(k+1)}$, and $S_{m_{k+1}}^{(k+1)}$ is the extra component. 
Let 
$C_u^{(k+1)}$ be the component which is directly outside of $S_{m_{k+1}}^{(k+1)}$. 
Let 
$\widetilde{S_d^{(k+1)}}$ ($\widetilde{S_u^{(k+1)}}$ resp.) be the union of the components of $\partial C_{m_{k+1}}^{(k+1)} \setminus S_{m_{k+1}}^{(k+1)}$ which are flat face down (up resp.) bowl like 2-spheres. 
Then $L^{(k)}$ is obtained from $L^{(k+1)}$ as follows. 
Let 
$\widetilde{B_d^{(k+1)}}$ ($\widetilde{B_u^{(k+1)}}$ resp.) be the union of mutually disjoint 3-ball(s) bounded by $\widetilde{S_d^{(k+1)}}$ ($\widetilde{S_u^{(k+1)}}$ resp.). 
First we remove $L^{(k+1)} \cap (C_u^{(k+1)} \cup C_{m_{k+1}}^{(k+1)} \cup \widetilde{B_d^{(k+1)}} \cup \widetilde{B_u^{(k+1)}})$ from $L^{(k+1)}$. 
Then add the pieces $L^{(k+1)} \cap \widetilde{B_d^{(k+1)}}$, $L^{(k+1)} \cap \widetilde{B_u^{(k+1)}}$, and $L^{(k+1)} \cap R_{m_{k+1}}^{(k+1)}$ by applying similarity deformations and parallel translations into a position appropriate for retrieving $L^{(k)}$. 
Finally add some monotonic arcs to obtain $L^{(k)}$. 
This observation implies the following.

\noindent
{\bf 
Facts 1.2B
}
\begin{enumerate}
\item 
If $i \ne u, m_{k+1}$, then the pairs $(R_i^{(k)}, L^{(k)} \cap R_i^{(k)})$ and $(R_i^{(k+1)}, L^{(k+1)} \cap R_i^{(k+1)})$ are similar. 

\item 
Note that $L^{(k+1)} \cap R_u^{(k+1)}$ consists of monotonic arcs. 
This fact allows us to regard $(R_u^{(k)}, L^{(k)} \cap R_u^{(k)})$ as obtained from $R_u^{(k+1)}$ and $R_{m_{k+1}}^{(k+1)}$ by putting them in a vertically disjoint position ($R_u^{(k+1)}$ is below, and $R_{m_{k+1}}^{(k+1)}$ is above) and adding some monotonic arcs. 
\end{enumerate}

See \figref{3.12}.

\begin{figure}[ht]\small\anchor{3.12}
\begin{center}
\psfraga <-4pt,-3pt> {R1}{$R_{m_{k+1}}^{(k+1)}$}
\psfraga <1pt,-5pt> {R2}{$R_u^{(k+1)}$}
\includegraphics[width=6cm, clip]{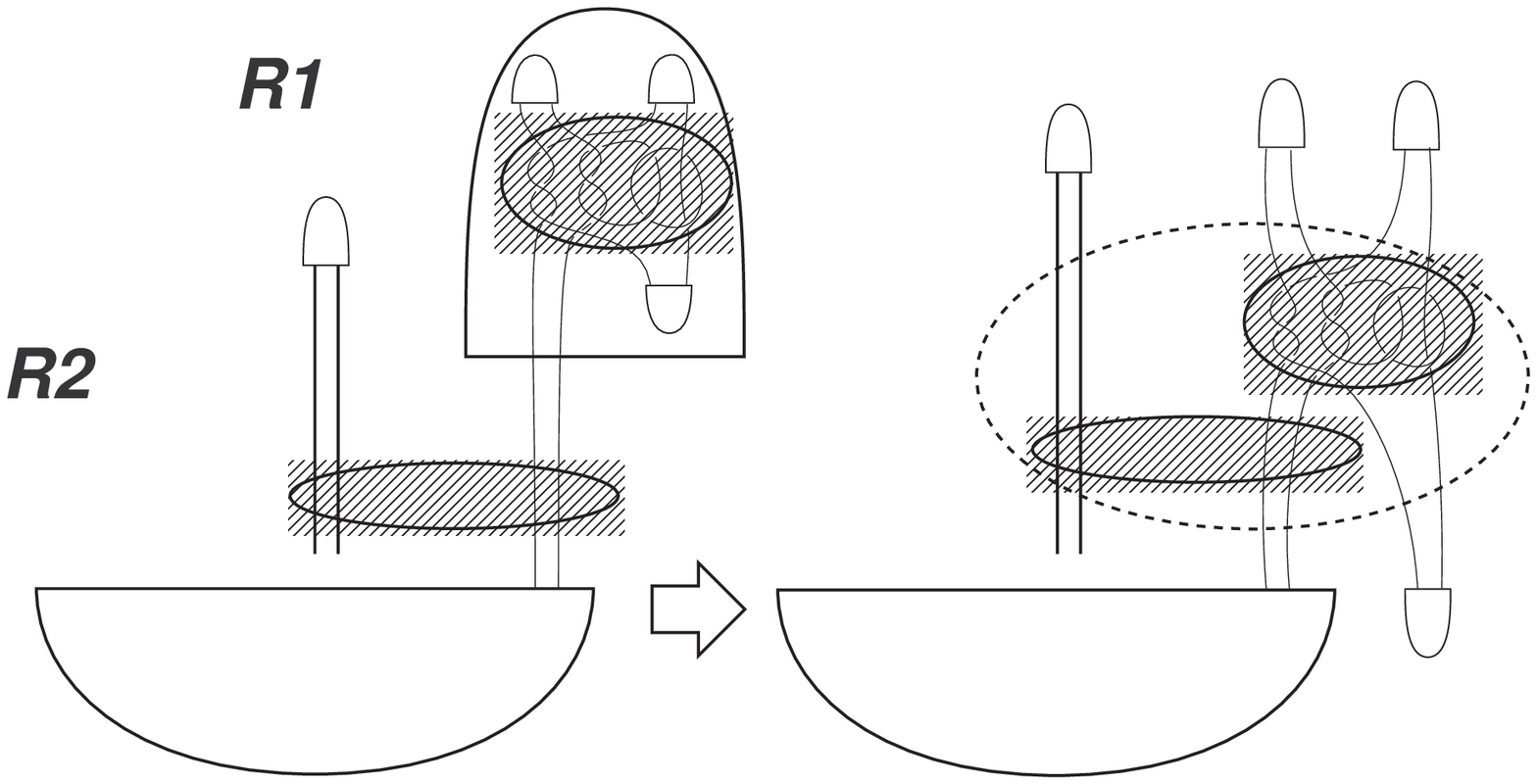}
\end{center}
\begin{center}
Figure 3.12
\end{center}
\end{figure}

For the remaining cases (Cases~2.1A$\sim$2.2B), analogous arguments apply, and detailed descriptions for the cases are left to the reader. 

\medskip
Suppose that $L$ has already been deformed as in Assertion~1. 
Let $f_t$ $(0 \le t \le 1)$ be the ambient isotopy from $L=L^{(0)}$ to $L^{(n)}$. 
Recall that $R_0^{(n)}, R_1^{(n)}, \dots , R_{m_n}^{(n)}$ are the cocoons of $L^{(n)}$ associated to $S^{(n)}$. 
Assertion~1 together with Facts~1.1A $\sim$ 2.2B implies the following.

\noindent 
{\bf 
Assertion 2}\qua
We have the following. 

\begin{enumerate}
\item 
For each $j$ $(j=0, 1, \dots , m_n)$, 
$L^{(n)} \cap R_j^{(n)}$ and 
$f_1^{-1}(L^{(n)} \cap R_j^{(n)})$ are similar. 

\item
The intervals 
$h(f_1^{-1}(L^{(n)} \cap R_0^{(n)}))$, 
$h(f_1^{-1}(L^{(n)} \cap R_1^{(n)}))$, \dots , 
$h(f_1^{-1}(L^{(n)} \cap R_{m_n}^{(n)}))$ 
are mutually disjoint. 
\end{enumerate}

By (2) of Assertion~2, we obtain a linear order, say $\prec$, 
on the cocoons $\{ R_0^{(n)}\!, R_1^{(n)}\!,$ $\dots , R_{m_n}^{(n)} \}$ that agrees with the positions of the intervals 
$h(f_1^{-1}(L^{(n)} \cap R_0^{(n)})),$ 
$h(f_1^{-1}(L^{(n)} \cap R_1^{(n)})), \dots , h(f_1^{-1}(L^{(n)} \cap R_{m_n}^{(n)}))$. 

By the description of the isotopies from $L^{(n)}$ to $L^{(0)}$ in Facts~1.1A$\sim$2.2B, we have the following.

\noindent
{\bf 
Assertion 3}\qua
The order $\prec$ is compatible with relative positions in $L^{(n)}$.

Let 
$G_0^{(n)}, G_1^{(n)}, \dots , G_{m_n}^{(n)}$ be the signed vertex graphs associated to $S^{(n)}$.

\noindent
{\bf 
Assertion 4}\qua
For each $G^{(n)}_j$, we have the following.

\begin{enumerate}
\item 
Each $G_j^{(n)}\backslash\{ \text{vertices} \}$ has either 

\begin{enumerate}
\item 
both a maximum and a minimum, or 
\item 
no critical points. 
\end{enumerate}

\item
$G_j^{(n)}$ is a bridge presentation giving a thin position. 
\end{enumerate}

\begin{proof}
Recall that $L^{(n)}$, $S^{(n)}$ satisfied Property~1 in Section~2 (see Section~3.2). 
The conclusion 1 is clear from Property~1 and the descriptions from the deformations described in 3.1. 
By (2), (3) of Property~1, we see that $G_j^{(n)}$ is in a bridge position. 
Suppose that the position $G_j^{(n)}$ is not thin. 
Let $G_j^{(n)}{}'$ be a thin position of the signed vertex graph $G_j^{(n)}$. 
(Hence $w(G_j^{(n)}{}') < w(G_j^{(n)})$.) 
Let $L^{(n)}{}'$ be the link obtained from $L^{(n)}$ by substituting $L^{(n)} \cap C_j^{(n)}$ with the 1-manifold corresponding to $G_j^{(n)}{}'$. 
Since $G_j^{(n)}{}'$ is isotopic to $G_j^{(n)}$, we may suppose that $L^{(n)}{}'$ is isotopic to $L^{(n)}$. 
Then consider the tangle decomposition of $L^{(n)}{}'$ by $S^{(n)}$. 
We may suppose that the cocoons of $L^{(n)}{}'$ are obtained from the cocoons of $L^{(n)}$ by replacing the cocoon $R_j^{(n)}$ with a cocoon, say $R_j^{(n)}{}'$, corresponding to $G_j^{(n)}{}'$. 
Recall that $\prec$ is the order of the cocoons $\{ R_0^{(n)}, R_1^{(n)}, \dots , R_{m_n}^{(n)} \}$ induced from the thin positon of $L$. 
By substituting $R_j^{(n)}$ to $R_j^{(n)}{}'$, we obtain an order $\prec'$ on the cocoons associated to $L^{(n)}{}'$, $S^{(n)}$. 
Clearly $\prec'$ is compatible with relative positions in $L^{(n)}{}'$. 
Then we apply the argument of the proof of Proposition~1 to obtain a position of the link $L$, say $L'$, which realizes the order. 
Note that $L'$ is obviously thinner than $L$, a contradiction. 
\end{proof}

\section{Main result}

In this section, we describe a search method for a thin position of a given non-splittable link $L$ satisfying the following (Assumptions~1, 2, 3). 

\noindent
{\bf 
Assumption 1}\qua
We know the bridge index, $n$, of $L$. 

\noindent
{\bf 
Assumption 2}\qua
We can give the list of all meridional, essential mutually non-parallel planar surfaces in $E(L)$ such that each planar surface has at most $2n-2$ boundary components. 

\noindent
{\bf 
Operation}\qua
Let $S=\cup_{i=1}^{m}S_i$ be a union of 2-spheres in $S^3$ such that $S \cap E(L)$ is a meridional, planar surface as in Assumption~2. 
Then, for each $i$ $(i=1, \dots , m)$, we assign $+$ to one side of $S_i$, and $-$ to the other. 
Note that there are $2^m$ ways to make such assignments. 
Let $C_0, C_1, \dots , C_m$ be the closures of the components of $S^3 \setminus \cup_{i=1}^m S_i$. 
Then for each $j$ $(j=0,1, \dots , m)$ the collar of each component of $\partial C_j$ has $+$ or $-$ sign. 
Then by regarding each component of $\partial C_j$ as a very tiny 2-sphere, we obtain a signed vertex graph, say $G_j$, from $L \cap C_j$. 
The third assumption is as follows. 

\noindent
{\bf 
Assumption~3}\qua
We know the bridge indices of all the signed vertex graphs obtained as above.

We now describe a method for obtaining various positions of $L$.

Recall that $n$ is the bridge index of $L$. 
Then an $n$-bridge presentation of $L$ is a candidate of a thin position of $L$. 
Here we note that the $n$-bridge presentation may be a thin position, 
even if $E(L)$ admits essential planar meridional surfaces (see Example~5.1 of \cite{H-K}). 
Let $S$ be a union of 2-spheres in $S^3$ such that $S \cap E(L)$ is a meridional, planar surface as in Assumption~2. 
Then we can obtain a number of systems of signed vertex graphs by using the procedures described in the operation above. 
Then, for each system of signed vertex graphs, we take minimal bridge presentations, say $G_0, G_1, \dots , G_m$, of the signed vertex graphs (Assumption~3). 
We expand the vertices of $G_0, G_1, \dots , G_m$ to make $+$ vertices ($-$ vertices resp.) flat face down (up resp.) bowl like 2-spheres. 
Then we combine the pieces, applying the inverse of deformations such as in \figref{2.3} to obtain a position of $L$, say $L'$, and a union of bowl like 2-spheres $S'$ with respect to which $L'$ satisfies Property~1 of Section 2.
Let $R_0, R_1, \dots , R_m$ be the cocoons of $L'$ associated to $S'$. 
Then consider all possible orders on $\{ R_0, R_1, \dots , R_m \}$ which are compatible with relative positions in $L'$. 
By Proposition~1 in Section 2, all such orders are realized as a position of $L$. 

Let ${\mathcal L}$ be the set of all positions obtained as above. 
Then the main result of this paper is as follows, 

\begin{theorem}
There is a thin position of $L$ in ${\mathcal L}$. 
\end{theorem}

\begin{proof}
Let $L^*$ be a thin position of $L$, and $\widetilde{L}$, $\widetilde{S}$ a position of $L$, and a union of bowl like 2-spheres giving an essential tangle decomposition of $\widetilde{L}$ obtained from $L^*$ by applying the procedure described in 3.2.

Let $\widetilde{S}=\cup\widetilde{S_i}$.  If $\widetilde{S_i}$ is flat face up, then we assign a + ($-$ resp.) to the inside collar (outside collar) of $\widetilde{S_i}$.  Similarly, we assign a $-$ (+ resp.) to the inside collar (outside collar) of $\widetilde{S_i}$ if it is flat face down.
Note that the union of planar surfaces $\widetilde{S} \cap E( \widetilde{L} )$, 
and the assignments of $\pm$ for the sides of $\widetilde{S}$ are included in one of the assignments of $\pm$ to the sides of planar surfaces described above. 
By Assertion~4 of Section 3.3, we see that the signed vertex graphs obtained from $\widetilde{L}$, $\widetilde{S}$ are in minimal bridge presentations. 
Possibly some of them are different from the bridge presentations introduced by Assumption~3. 
However the corresponding positions must have the same bridge index, because they are minimal bridge presentations of the same signed vertex graph. 
Let $\widetilde{L}'$ be the position of $L$ obtained by combining the minimal bridge presentations obtained by Assumption~3. 
It is clear that the relative positions of the cocoons of $\widetilde{L}$ associated to $\widetilde{S}$ coincides with the relative positions of the cocoons of $\widetilde{L}'$ associated to $\widetilde{S}$. 
Note that the ordering of the cocoons for $\widetilde{L}$ coming from $L^*$ is compatible with the relative positions of $\widetilde{L}$ by Facts~1.1A$\sim$2.2B of Section 3.3.
One ordering of cocoons for $\widetilde{L}'$ is exactly the ordering of $\widetilde{L}$ coming from $L^*$, and it is compatible with relative positions.
By applying Proposition~1 of Section 2 to $\widetilde{L}'$, we obtain a position of $L$, say $\widetilde{L}^*$, 
which has the same width data as $L^*$, hence $\widetilde{L}^*$ is in a thin 
position. 
This completes the proof of the theorem. 
\end{proof}

\begin{remark}
As a practical application of the above procedures, we can reduce the number of candidates of thin positions by using various observations (e.g. Assertion~4). 
For example, see Section~5 of this paper. 
\end{remark}

\section{Example}

In this section, by using Theorem~2, we show that the thin position for the pretzel link $L = P(3,3,3,3,3,3)$ is essentially as in \figref{5.1} (2), whose width is 48. 

\begin{figure}[ht]\small\anchor{5.1}
\begin{center}
\psfrag{L}{$L$}
\includegraphics[width=6cm, clip]{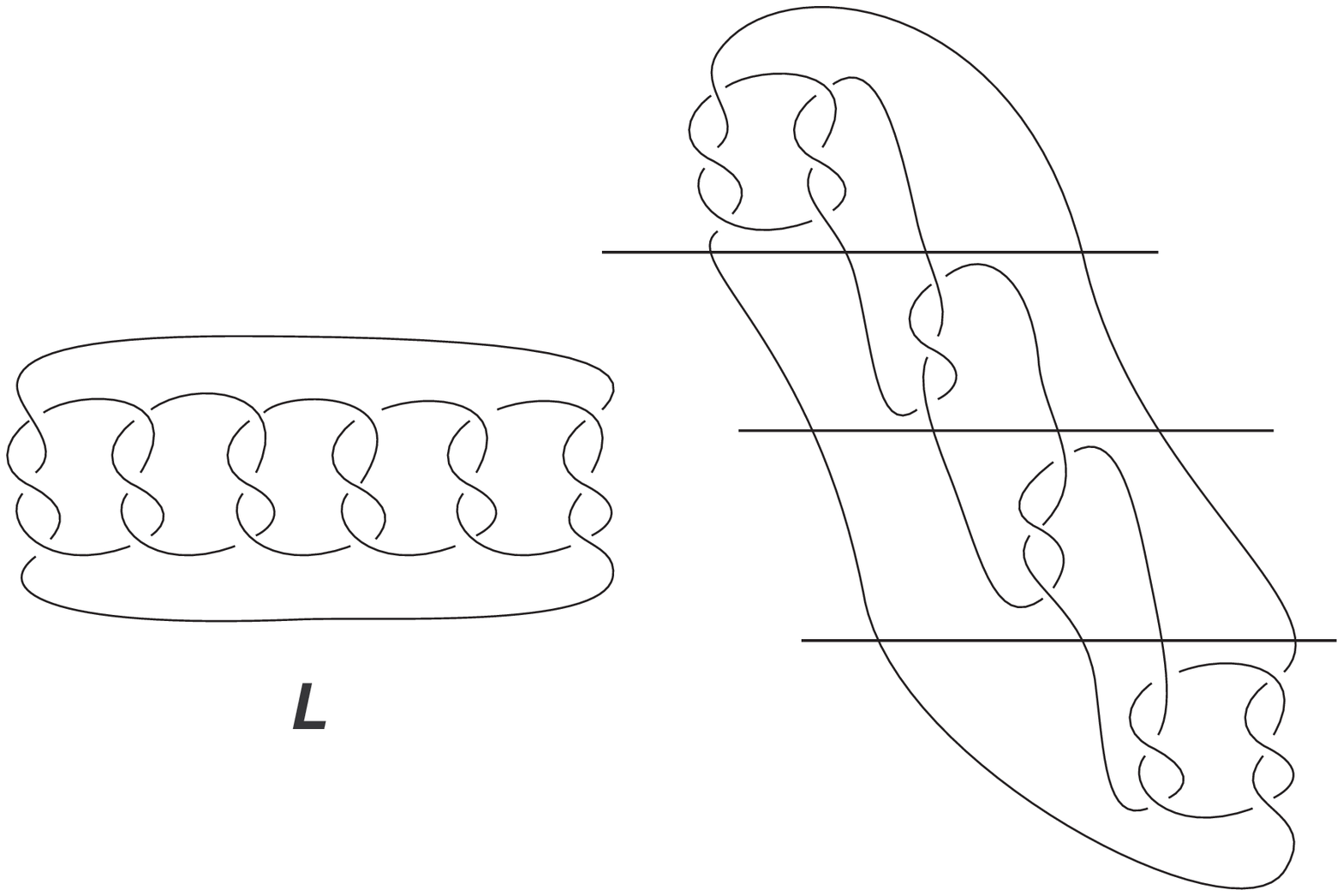}
\end{center}
\begin{center}
Figure 5.1
\end{center}
\end{figure}

Let $M$ be the 2-fold branched covering space of $S^3$ with branch set $L$. 
We note that $M$ is a Seifert fibered space (see, for example, 12.30 of \cite{B-Z}). 
In particular, $M=\{0,0,(3,1),(3,1),(3,1),(3,1),(3,1),(3,1)\}$.

\noindent
{\bf Claim 1}\qua 
The bridge index of $L$ is $6$.

\begin{proof}
It is elementary to check, using \cite{M-S}, that $M$ has no horizontal Heegaard splittings, and thus that its Heegaard genus is 5. 
Then, by 11.5 of \cite{B-Z}, we see that the bridge index for $L$ is at least 6. 
It is directly observable from \figref{5.1} that $L$ actually admits a 6 bridge presentation. 
Hence the bridge index of $L$ is 6. 
\end{proof}

\noindent
{\bf Claim 2}\qua 
Suppose $L$ is in a thin position. 
Then there does not exist a thin 2-sphere $S$ intersecting $L$ in more than four points.

\begin{proof}
We first show: 

\noindent
{\bf Subclaim 2.1}\qua 
In general, if there exist inequivalent thick 2-spheres $S_1$, $S_2$ for $L$, and $S_i$ intersects $L$ in $2n_i$ points, then the width of $L$ is greater than or equal to $n_1(n_1 + 1) + n_2(n_2 + 1)$.

\begin{proof}
Without loss of generality, we may assume that $S_1$ is higher than $S_2$. 
Then there is a sequence of (perhaps not successive) 2-spheres $S_1^{(1)}, \dots , S_1^{(n_1 - 1)}$ above $S_1$ intersecting $L$ transversely in $2,$ $4,$ $\dots$, $2(n_1 - 1)$ points respectively. 
Hence these 2-spheres together with $S_1$ contribute $2 +$$4 +$$\dots + 2n_1$ $= n_1 (n_1+1)$. 
Apply the same argument to the region below $S_2$. 
These estimations show that $w(L) \ge n_1(n_1 + 1) + n_2(n_2 + 1)$. 
\end{proof}

Assume that there exists a thin 2-sphere $S$ which intersects $L$ in $2n$ points with $n \ge 4$. 
Then there exist thick 2-spheres $S_1$, $S_2$ such that $S_1$ ($S_2$ resp.) is above (below resp.) $S$ and $S_i$ $(i=1,2)$ intersects $L$ in $2 n_i$ points with $n_i > n$. 
By Subclaim~2.1, we see that the width of $L$ is greater than or equal to $2 \times (n+1)(n+2)$ which is greater than or equal to $2 \times (4+1)(4+2) = 60$. 
However, by \figref{5.1}, we see that $L$ admits a presentation with width 48, a contradiction. 

Hence we may assume that $n=3$. 
Then the least width we can construct from this thin 2-sphere is $2 + 4 + 6 + 8 + 6 + 8 + 6 + 4 + 2 = 46$. 
However this would imply that the presentation of $L$ has exactly 5 maxima (and minima)(\figref{5.2}), contradicting Claim~1. 
Any other width constructed from this thin 2-sphere is greater than 48, a contradiction. \end{proof}

\begin{figure}[ht]\small\anchor{5.2}
\begin{center}
\psfrag{8-braid}{8-braid}
\includegraphics[width=3cm, clip]{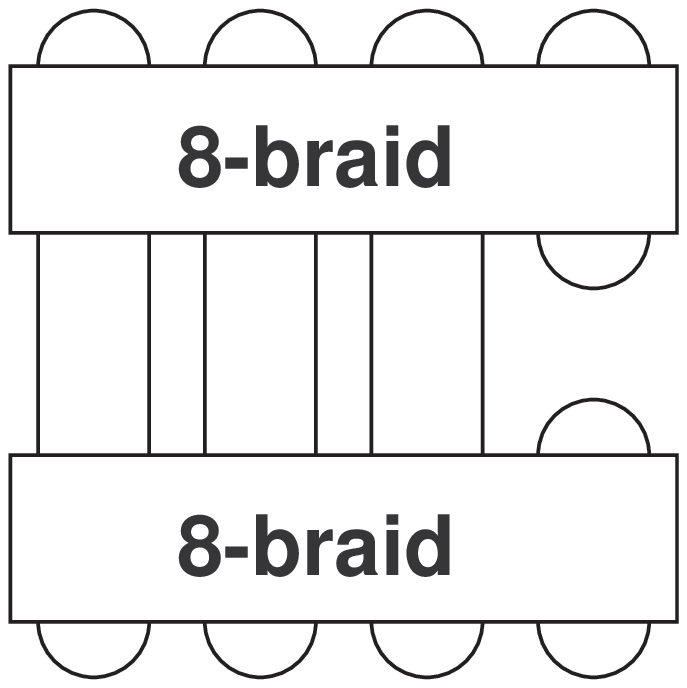}
\end{center}
\begin{center}
Figure 5.2
\end{center}
\end{figure}

Since $L$ is not composite, $n\ne1$.
Thus for the purpose of searching for a thin position, it is enough to consider only the case of $n=2$.

We assume that the reader is familiar with \cite{O}. 
By using the presentation of the fundamental group of $M$ in 11.31 of \cite{B-Z}, we see that $M$ does not have positive Betti number. 
This implies that $L$ does not admit a surface that is horizontal with respect to the orbifold fibration structure on $(S^3,L)$.
Then by 2.8 and 2.14 of \cite{O}, we see that each incompressible, $\partial$-incompressible, meridional, planar surface in the exterior of $L$ with at most four boundary components is one of \figref{5.3}, up to homeomorphisms of $E(L)$. 

\begin{figure}[ht]\small\anchor{5.3}
\begin{center}
\includegraphics[width=7cm, clip]{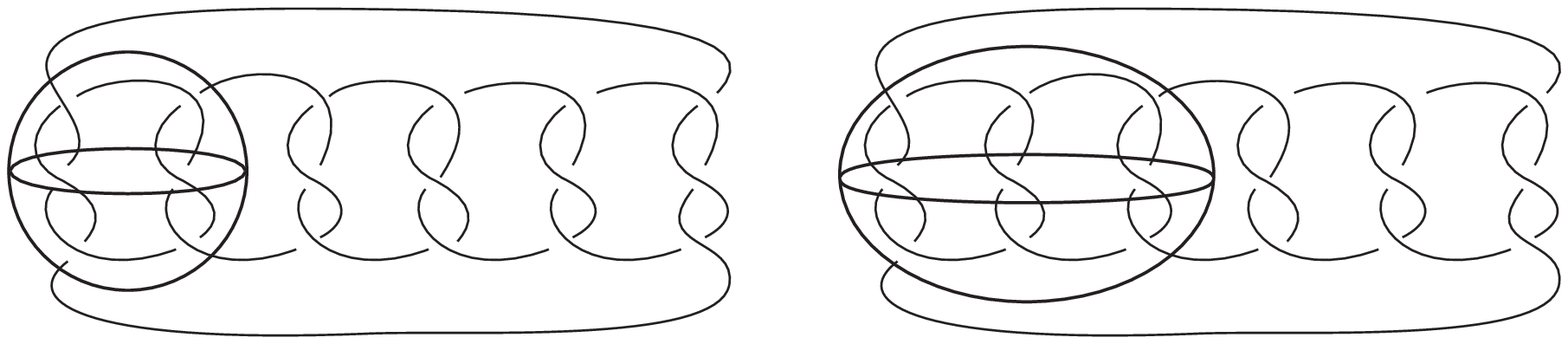}
\end{center}
\begin{center}
Figure 5.3
\end{center}
\end{figure}

\begin{figure}[b!]\small\anchor{5.4}
\begin{center}
\psfrag{1}{1}
\psfrag{2}{2}
\psfrag{3}{3}
\psfrag{4}{4\qquad\qquad5\qquad\qquad6\qquad\qquad7}
\psfrag{8}{8\qquad\qquad\qquad\ \ \ 9}
\psfrag{10}{10\qquad\qquad11}
\includegraphics[width=7cm, clip]{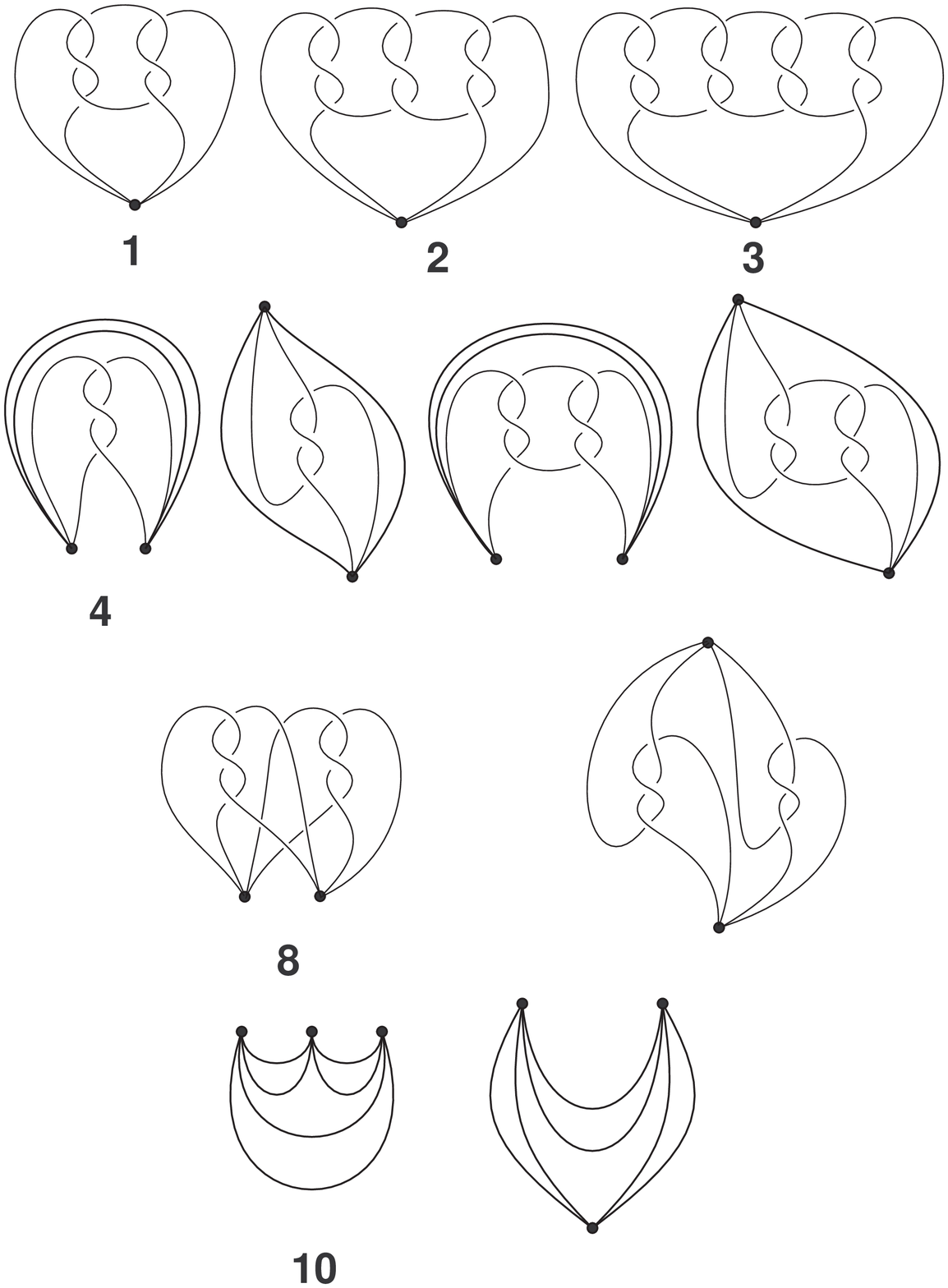}
\end{center}
\begin{center}
Figure 5.4
\end{center}
\end{figure}

Then each spacial graph obtained from $L$ as in Section~4 is one of \figref{5.4}. 
It is easy to see that each graph equipped with sign has (minimal) bridge presentations as in \figref{5.4}, since if not, any smaller bridge presentation would induce a smaller bridge presentation for $L$.
By Assertion~4 of Section~3.3, we see that we can throw graphs 4, 8, 10 and 11 of \figref{5.4} out of consideration. 
It is easy to see that the signed vertex graphs 2, 3, 7, 9 admit thinner presentation 
(which is not a bridge presentation). 
Hence we may throw them out of consideration by Assertion~4 of Section~3. 
Hence we can apply the arguments of Section~4 to the graphs 1, 5, and 6, and this gives the thin position of $L$ is as in \figref{5.1}.

\Addresses\recd

\end{document}

%% file: agtout.tex
\def\ifplaintex{\expandafter\ifx\csname documentclass\endcsname\relax}

\def\gtp{{\mathsurround=0pt\it $\cal G\mskip-2mu$eometry \&\ 
$\cal T\!\!$opology $\cal P\!$ublications}}  

\def\recd{{\small Received:\qua\receiveddate\ifx\reviseddate\relax
\else\qquad Revised:\qua\reviseddate\fi\par}} 

\def\lognumber#1{\def\thelognumber{#1}}
\def\volumenumber#1{\def\thevolumenumber{#1}}
\def\volumeyear#1{\def\thevolumeyear{#1}}
\def\papernumber#1{\def\thepapernumber{#1}}
\def\pagenumbers#1#2{\def\startpage{#1}\def\finishpage{#2}}
\def\published#1{\def\publishdate{#1}}

\def\received#1{\def\receiveddate{#1}}

\def\accepted#1{\def\accepteddate{#1}}

\def\asciiaddress#1{\def\theasciiaddress{#1}}
\def\asciiemail#1{\def\theasciiemail{#1}}

\let\\\par\let\thelognumber\relax\let\thevolumenumber\relax
\let\thepapernumber\relax\let\thevolumeyear\relax\let\startpage\relax
\let\finishpage\relax\let\publishdate\relax\let\receiveddate\relax
\let\reviseddate\relax\let\accepteddate\relax\let\theasciititle\relax
\let\theasciiauthors\relax\let\theasciiaddress\relax
\let\theasciiabstract\relax

\let\theasciiemail\relax

\ifplaintex
\font\logobig=cmssbx10 scaled 3836
\font\logomed=cmssbx10 scaled 2557
\else
\font\logobig=cmssbx10 scaled 4200
\font\logomed=cmssbx10 scaled 2800
\fi

\long\def\makeagttitle{   
\count0=\startpage
\agt\hfill      
\hbox to 45truept{\vbox to 0pt{\vglue -13truept{\logomed A\kern -.37em{\logobig 
T}\kern -.38em G}\vss}\hss}
\break
{\small Volume \thevolumenumber\ (\thevolumeyear)
\startpage--\finishpage\nl
Published: \publishdate}

\vglue .25truein

{\parskip=0pt\leftskip 0pt plus
1fil\def\\{\par\smallskip}{\Large\bf\thetitle}\par\medskip} \vglue
0.05truein

{\parskip=0pt\leftskip 0pt plus 1fil\def\\{\par}{\sc\theauthors}
\par\medskip}%
 
\vglue 0.03truein 

{\small\leftskip 25truept\rightskip 25truept{\bf Abstract}\stdspace\theabstract

{\bf AMS Classification}\stdspace\theprimaryclass
\ifx\thesecondaryclass\relax\else; \thesecondaryclass\fi\par
{\bf Keywords}\stdspace \thekeywords\par}\vglue 7truept

}   

\ifplaintex
\font\phead=cmsl9 scaled 950
\font\pnum=cmbx10 scaled 913
\font\pfoot=cmsl9 scaled 950
\headline{\vbox to 0pt{\vskip -4.5mm\line{\small\phead\ifnum
\count0=\startpage ISSN 1472-2739 (on-line) 1472-2747 (printed)
\hfill {\pnum\folio}\else\ifodd\count0\def\\{ }%
\ifx\theshorttitle\relax\thetitle\else\theshorttitle\fi\hfill{\pnum\folio}
\else\def\\{ and }{\pnum\folio}\hfill\ifx\theshortauthors\relax\theauthors
\else\theshortauthors\fi\fi\fi}\vss}}
\footline{\vbox to 0pt{\vglue 0mm\line{\small\pfoot\ifnum\count0=\startpage
\copyright\ \gtp\hfill\else
\agt, Volume \thevolumenumber\ (\thevolumeyear)\hfill\fi}\vss}}
\else
\font=cmsl9 scaled 1050
\font\lnum=cmbx10 
\font=cmsl9 scaled 1050
\makeatletter
\def\@oddhead{{\small\lhead\ifnum\count0=\startpage ISSN 1472-2739 
(on-line) 1472-2747 (printed)\hfill {\lnum\number\count0}\else\ifodd\count0
\def\\{ }\ifx\theshorttitle\relax \thetitle \else\theshorttitle\fi\hfill
{\lnum\number\count0}\else\def\\{ and }{\lnum\number\count0}
\hfill\ifx\theshortauthors\relax 
\theauthors\else\theshortauthors\fi\fi\fi}}\def\@evenhead{\@oddhead}
\def\@oddfoot{\small\lfoot\ifnum\count0=\startpage\copyright\ \gtp\hfill\else
\agt, Volume \thevolumenumber\ (\thevolumeyear)\hfill\fi}
\def\@evenfoot{\@oddfoot}
\makeatother
\fi
\let\maketitlepage\makeagttitle
\let\makeshorttitle\maketitlepage
\let\maketitle\maketitlepage

\newwrite\gtoutfile
\long\gdef\makeheadfile{  
{\def\\{, }\def\s{ }
\immediate\openout\gtoutfile head.xxx
\immediate\write\gtoutfile{Proxy-for: \ifx\theasciiauthors\relax
\theauthors\else\theasciiauthors\fi\s<\ifx\theasciiemail\relax\theemail\else\theasciiemail\fi>}
\immediate\write\gtoutfile{\noexpand\\}
\immediate\write\gtoutfile{Authors: \ifx\theasciiauthors\relax
\theauthors\else\theasciiauthors\fi}
{\def\\{ }\immediate\write\gtoutfile{Title: \ifx\theasciititle\relax
\thetitle\else\theasciititle\fi}}
\immediate\write\gtoutfile{Subj-class: GT or SG, GR etc}
\immediate\write\gtoutfile{MSC-class: \theprimaryclass\ifx\thesecondaryclass\relax\else, \thesecondaryclass\fi}
\immediate\write\gtoutfile{Journal-ref: Algebr. Geom. Topol. \thevolumenumber\s
(\thevolumeyear) \startpage-\finishpage}
\immediate\write\gtoutfile{Comments: Published by Algebraic and
Geometric Topology at}
\immediate\write\gtoutfile{\s\s\s  http://www.maths.warwick.ac.uk/agt/AGTVol\thevolumenumber/agt-\thevolumenumber-\thepapernumber.abs.html}
\immediate\write\gtoutfile{\noexpand\\}
\immediate\write\gtoutfile{}
\ifx\theasciiabstract\relax
\immediate\write\gtoutfile{\theabstract}\else
\immediate\write\gtoutfile{\theasciiabstract}\fi
\immediate\write\gtoutfile{}
\immediate\write\gtoutfile{\noexpand\\}
\immediate\write\gtoutfile{}
\immediate\closeout\gtoutfile}}  

\def\maketitlepage{\makeagttitle\makeheadfile}
\let\makeshorttitle\maketitlepage
\let\maketitle\maketitlepage